\documentclass[12pt]{amsart}
\usepackage{fullpage,enumerate}
\usepackage{amsmath}
\usepackage{graphicx,amssymb,amsthm,tikz}
\usepackage{hyperref}
\usepackage[font=footnotesize,labelfont=bf]{caption}
\usepackage{tikz}
\usepackage{multirow}

\newcommand{\affil}[1]{}

\allowdisplaybreaks[3] 
\everymath{\displaystyle} 

\newcommand{\excep}{exceptional pattern}
\newcommand{\exceps}{exceptional patterns}

\newcommand{\ASM}{\textsf{\textit{ASM}}}

\newtheorem{thm}{Theorem}[section]
\newtheorem{prop}[thm]{Proposition}
\newtheorem{lem}[thm]{Lemma}

\newtheorem*{MainTheorem}{Theorem~\ref{thm:superexponential}}
\newtheorem*{BetterBoundThm}{Theorem~\ref{T:3412}}
\theoremstyle{definition}
\newtheorem*{ThreeNegativeOnesThm}{Theorem~\ref{thm:exactly3enumeration}}
\newtheorem*{kNegativeOnesThm}{Theorem~\ref{thm:knegativeones}}

\theoremstyle{definition}

\newtheorem{remark}[thm]{Remark}

\numberwithin{equation}{section}

\title[Enumeration of pattern-avoiding alternating sign matrices: An asymptotic dichotomy]{Enumeration of pattern-avoiding alternating sign matrices: An asymptotic dichotomy}

\author{Mathilde Bouvel\affil{1}}
\address{\affil{1} Universit\'{e} de Lorraine \\
CNRS, Inria, LORIA, F-54000 \\
Nancy, France}
\email{mathilde.bouvel@loria.fr}

\author{Eric S. Egge\affil{2}}
\address{\affil{2} Department of Mathematics and Statistics \\
Carleton College \\
Northfield, MN, USA}
\email{eegge@carleton.edu}

\author{Rebecca N. Smith\affil{3}}
\address{\affil{3} Department of Mathematics \\
SUNY Brockport \\
Brockport, NY, USA}
\email{rnsmith@brockport.edu}

\author{Jessica Striker\affil{4}}
\address{\affil{4} Department of Mathematics \\
North Dakota State University \\
Fargo, ND, USA}
\email{jessica.striker@ndsu.edu}

\author{Justin M. Troyka\affil{5}}
\address{\affil{5} Department of Mathematics \\
California State University, Los Angeles \\
CA, USA}
\email{jtroyka@calstatela.edu}

\date{September 12, 2025}

\begin{document}

\begin{abstract}We completely classify the asymptotic behavior of the number of alternating sign matrices classically avoiding a single permutation pattern, in the sense of [Johansson and Linusson 2007]. In particular, we give a uniform proof of an exponential upper bound for
the number of alternating sign matrices classically avoiding one of eleven particular patterns,
and a super-exponential lower bound for all other single-pattern avoidance classes. 
We also show that for any fixed integer $k$, there is an exponential upper bound for the number of alternating sign matrices that classically avoid any single permutation pattern and contain precisely $k$ negative ones. Finally, we prove that there must be at most $3$ negative ones in an alternating sign matrix which classically avoids both $2143$ and $3412$, and we exactly enumerate the number of them with precisely $3$ negative ones. 
\end{abstract}

\maketitle
\section{Introduction}
Permutations are fundamental objects in combinatorics and algebra. The study of various patterns in permutations and the set of permutations that avoid certain patterns has been a fruitful area of research (see e.g.~\cite{Bevan,Kitaev,Vatter} \cite[Ch.\ 4]{Bona}). Alternating sign matrices are an interesting superset of permutations with connections to algebra, geometry, dynamics, and statistical physics (see e.g.~\cite{BressoudBook,FK_ASM,Striker}).   This paper studies an analog of pattern avoidance on alternating sign matrices, proving a dichotomy in the growth rates of their pattern avoidance classes. We begin with some definitions, which we discuss in more detail in Section~\ref{sec:def}.

An {\em alternating sign matrix} (ASM) is a square matrix such that
\begin{itemize}
\item
every entry is $0$, $1$, or $-1$,
\item 
in each row and in each column the sum of the entries is $1$, and
\item 
in each row and in each column the nonzero entries alternate in sign.
\end{itemize}
In this context, a {\em permutation matrix} (or permutation) is an ASM in which every entry is $0$ or $1$.

If $A$ is an $n \times n$ ASM and $B$ is a $k \times k$ permutation, then we say $A$ {\em contains} $B$ whenever there is a $k \times k$ submatrix $C$ of $A$ such that for all $i$ and $j$, if $B_{ij} = 1$ then $C_{ij} = 1$.
In this case we sometimes call $C$ an {\em occurrence} of $B$ in $A$.
We note, however, that $C$ can be an occurrence of $B$ without being an exact copy of $B$, since any entry of $C$ whose corresponding entry of $B$ is $0$ can be $0$, $1$, or $-1$.
For example, if $A$ and $B$ are given by
\[ A= \left( \begin{matrix} 0 & 0 & 0 & 1 & 0 & 0 \\ 0 & {\bf 1} & {\bf 0} & {\bf 0} & 0 & {\bf 0} \\ 1 & -1 & 1 & -1 & 1 & 0 \\ 0 & {\bf 1} & {\bf -1} & {\bf 1} & -1 & {\bf 1} \\ 0 & {\bf 0} & {\bf 1} & {\bf -1} & 1 & {\bf 0} \\ 0 & {\bf 0} & {\bf 0} & {\bf 1} & 0 & {\bf 0} \end{matrix} \right); \hspace{30pt} B=\left(\begin{matrix} 1 & 0 & 0 & 0 \\ 0 & 0 & 0 & 1 \\ 0 & 1 & 0 & 0 \\ 0 & 0 & 1 & 0\end{matrix}\right); \]
then rows $2,4,5,6$ and columns $2,3,4,6$ of $A$ are an occurrence of $B$ in $A$.
The entries in the resulting submatrix are in bold in $A$ above.

We are interested in the case in which an ASM $A$ does not contain a permutation $B$;  in this case we say $A$ {\em classically avoids} $B$.
This notion of avoidance in ASMs was first studied by Johansson and Linusson~\cite{JL}; we call it {\em classical avoidance} to distinguish it from the notion of {\em key-avoidance} studied in~\cite{keyAv}. This notion is also distinct from the \emph{zero-nonzero patterns} studied in \cite{BKMS}.
For a set $R$ of permutations, we write $\ASM_n(R)$ to denote the set of $n \times n$ ASMs classically avoiding all permutations in $R$, and $\ASM(R) = \bigsqcup_{n\ge0} \ASM_n(R)$.
The main results of~\cite{JL} were enumerations of the sets $\ASM_n(132)$ by the large Schr\"oder numbers and $\ASM_n(132,123)$ by every second Fibonacci number. 
In Table~\ref{tab:first_values}, we give the initial terms of the enumeration sequences for ASMs classically avoiding a single pattern of length three or four that were originally computed by Johansson and Linusson~\cite{JL}.  Note symmetries from rotations and reflections, which include permutation reversal, complementation, and inverse  yield identical results. 
In this paper, we primarily consider growth rates instead of exact enumerations (with the notable exception of Theorem~\ref{thm:exactly3enumeration}, which we describe at the end of the introduction).

\begin{table}[ht]
    \centering
    \begin{tabular}{|l|l|l|}
    \hline 
       Excluded pattern  
       & First terms of enumeration sequence  & Growth of the sequence,\\
        up to symmetry &   & exact enumeration \\
    \hline 
        \multirow{ 3}{*}{132} & \multirow{3}{*}{1, 2, 6, 22, 90, 394, 1806, 8558,\ldots} & exponential, specifically  \\
        &&  \href{https://oeis.org/A000108}{OEIS A000108} \emph{i.e.}\ large \\
        && Schr\"{o}der numbers, see~\cite{JL}\\  
    \hline 
        123 & 1, 2, 6, 23, 103, 514, 2785, 16132, \ldots & super-exponential, open\\
    \hline 
        2143 & 1, 2, 7, 40, 320, 3152, 35551, 441280, \ldots & exponential, open\\
    \hline 
        2413 & 1, 2, 7, 41, \textbf{364}, 4168, 54659, 775528, \ldots & exponential, open\\
    \hline 
        1243 & 1, 2, 7, 41, 360, 4200, 59869, 990930, \ldots & super-exponential, open\\
    \hline 
        1432 & 1, 2, 7, 41, 361, 4234, 60723, 1009328, \ldots & super-exponential, open\\
    \hline 
        1342 & 1, 2, 7, 41, 367, 4455, 66403, 1138774, \ldots & super-exponential, open\\
    \hline 
        1234 & 1, 2, 7, 41, 370, 4638, 74093, 1423231, \ldots & super-exponential, open\\
     \hline 
        1324 & 1, 2, 7, 41, 376, 4985, 88985, 2024954, \ldots & super-exponential, open\\
    \hline 
    \end{tabular}
    \caption{Enumeration sequences and growth rates (from Theorem~\ref{thm:superexponential}) of ASMs classically avoiding a permutation pattern. Every permutation pattern of size $3$ or $4$ is listed, up to symmetries of the square. As an aside, the bold entry in the $2413$ sequence is greater than the corresponding number for $1243$ and $1432$, even though $\ASM(2413)$ is exponential and the other two are super-exponential.
    }
    \label{tab:first_values}
\end{table}

Classical avoidance in ASMs reduces to the classical notion of pattern avoidance in permutations when $A$ is a permutation matrix, and in this context it is known that every set of pattern-avoiding permutations grows at most exponentially.
More specifically, we can extend our definition of permutation containment and avoidance to the set of all matrices of $0$s and $1$s.
In this context, Marcus and Tardos have shown that for any permutation $\sigma$, there is a constant $c(\sigma)$ such that the number of $n\times n$ $0-1$ matrices which avoid $\sigma$ is bounded above by $c(\sigma)^n$~\cite{MarcusTardos}.
For permutations more is true:  if $S_n(\sigma)$ is the set of $n \times n$ permutations which avoid $\sigma$, then $\lim_{n\to\infty} \sqrt[n]{|S_n(\sigma)|}$ exists~\cite{Arratia}.
We sometimes call this value the {\em Stanley--Wilf limit} of the class $S_n(\sigma)$, or simply the Stanley--Wilf limit of $\sigma$.
Table \ref{table:SWlimits} gives some of the few results that are known about the values of specific Stanley--Wilf limits.
More information on known Stanley--Wilf limits, including examples demonstrating that some Stanley--Wilf limits are not integers and not rational, can be found in \cite{Bonairrational} and \cite{Bonarecords}.
\begin{center}
    \begin{table}[h]
    \begin{tabular}[h]{c|c|c}
    $\sigma$ & Stanley--Wilf limit of $\sigma$ & References \\
    \hline
    123 & \multirow{2}{*}{4} & \multirow{2}{*}{\cite{Knuth,MacMahon, Regev,  SS85}}\\
    132 & \\
    \hline
    1342 & \multirow{2}{*}{8} & \multirow{2}{*}{\cite{BloomElizalde,Bona1342}} \\
    2413 & \\
    \hline
    1234 & \multirow{4}{*}{9} & \multirow{4}{*}{\cite{BWX,Bona,Regev}} \\
    1243 & \\
    1432 & \\
    2143 & \\
    \hline
    1324 & between 10.271 and 13.5 & \cite{Bevan1324} \\
    \hline
    $12\cdots k$ & $(k-1)^2$ & \cite{Regev}
    \end{tabular}
        \caption{What is known about some Stanley--Wilf limits.}
        \label{table:SWlimits}
    \end{table}
\end{center}

Returning to the more general case, we note that Johansson and Linusson's results imply that $\lim_{n \to \infty} \sqrt[n]{|\ASM_n(132)|} = 3+2\sqrt{2}$ and $\lim_{n \to \infty} \sqrt[n]{|\ASM_n(132,123)|} = \frac{3+\sqrt{5}}{2}$, so both of these classes of pattern-avoiding ASMs grow exponentially.
In this paper we find that this is unusual for pattern-avoiding ASMs, by showing first that $132$ is one of just eleven\footnote{We may add a twelfth permutation to these exceptions, namely the empty permutation. However, as avoiding the empty permutation is both uninteresting and confusing, we have chosen to restrict to excluded patterns of positive size.} exceptional permutations whose classical pattern avoidance class grows (at most) exponentially, and second that for all other permutations the corresponding class grows super-exponentially.
More specifically, let $X$ denote the following set of permutations:
\[ X := \{1, 12, 21, 132, 213, 231, 312, 2143, 2413, 3142, 3412\}.\]
We note that, up to symmetries (i.e.\@ regarding the reverse or inverse of a pattern as equivalent to the pattern), $X$ contains only five essentially different permutations. 
We also note that these permutations can be characterized as the permutations which avoid both $123$ and $321$, or as the permutations which can be constructed by choosing points on a rectangle whose sides are exactly vertical or horizontal (this is explained in further detail in the proof of Theorem~\ref{thm:superexponential}). 
Our first main theorem is the following.

\begin{thm} \label{thm:superexponential}
Suppose $\sigma$ is a permutation. 
If $\sigma\in X$, then there exists $c$ such that $|\ASM_n(\sigma)| \leq c^n$ for all $n$. Otherwise,
 $|\ASM_n(\sigma)| \ge \lfloor n/3 \rfloor !$ for all $n$.
\end{thm}

We note that this theorem gives two meanings to the notation $X$ above: these \textbf{ex}ceptional patterns are exactly those whose avoidance describes a family of ASMs with (at most) \textbf{ex}ponential growth. Table~\ref{tab:first_values} gives a visual summary of this theorem by including the information for each pattern of whether the growth rate is exponential or super-exponential.

In interesting contrast, we show that for any permutation $\sigma$, and for any fixed integer $k$, the sequence enumerating $n\times n$ ASMs that classically avoid $\sigma$ and contain exactly $k$ negative ones grows at most exponentially with $n$. Let  $\ASM_{n,k}(\sigma)$ denote the set of ASMs in $\ASM_n(\sigma)$ containing exactly $k$ negative ones.

\begin{thm} \label{thm:knegativeones}
    Suppose $\sigma$ is a permutation and $k$ is a non-negative integer.  
Then there exists $c$ such that $|\ASM_{n,k}(\sigma)| \leq c^n$ for all $n$. 
\end{thm}

Our next main result is an explicit upper bound, achieved by construction, for the class of ASMs avoiding a specific pattern in $X$. 

\begin{thm}~\label{T:3412}
The number of alternating sign matrices that classically avoid $3412$ has the following exponential upper bound:
$|\ASM_n(3412)| \leq 62208^n$. The same holds for  $\ASM_n(2143)$.
\end{thm}

Our final main theorem is an exact enumeration of ASMs with three negative ones avoiding both of the patterns $2143$ and $3412$. 
    \begin{thm}
    \label{thm:exactly3enumeration}
        For $n\geq 7$, we have
        \[ |\ASM_{n,3}(2143, 3412)| = \sum_{k=0}^{n-5} 2^{n-4-k} \binom{2k}{k} - (n-2) 2^{n-5}. \]
        For $n \le 6$, we have $|\ASM_{n,3}(2143, 3412)| = 0$. The generating function for this sequence is
        \[ \sum_{n\ge0} |\ASM_{n,3}(2143, 3412)|\,z^n =  \frac{2z^5}{(1-2z)\sqrt{1-4z}} -\frac{2z^5(1-2z^2)}{(1-2z)^2}. \]
    \end{thm}
This result is motivated by the fact, which we prove in Lemma~\ref{lem:atmost3}, that any ASM in $\ASM_n(2143,3412)$ has at most $3$ negative ones. 
In particular, it represents a first step toward an exact enumeration of $\ASM_n(2143,3412)$.
As a second step toward this enumeration, we note that the class $S_n(2143,3412)$, which is the set of permutations in $\ASM_n(2143,3412)$, is the class of so-called \emph{skew-merged} permutations, which were first enumerated by Atkinson \cite{Atkinson}.

The paper is organized as follows. Section~\ref{sec:def} specifies  definitions and notation. Section~\ref{sec:mainthm} proves Theorem~\ref{thm:superexponential}, generalizes it in Theorem~\ref{thm:incdec} to certain sets of patterns, and proves Theorem~\ref{thm:knegativeones}. Section~\ref{sec:abetterbound} proves Theorem~\ref{T:3412} and some useful lemmas. Section~\ref{sec:toward} studies the exact enumeration of $\ASM_n(2143,3412)$, proving Theorem~\ref{thm:exactly3enumeration}. We end in Section~\ref{sec:future} with some ideas for future work. 

\section{Definitions and notation}
\label{sec:def}
We let $S_n$ denote the set of permutations of size $n$, and we define $S = \bigsqcup_{n\ge0} S_n$. We identify each $\pi \in S_n$ with its permutation matrix
\[ P_\pi = \begin{bmatrix} e_{\pi(1)} \\ e_{\pi(2)} \\ \cdots \\ e_{\pi(n)} \end{bmatrix}, \]
where $e_i$ is the row vector with a $1$ in column $i$ and $0$'s everywhere else. For example, the permutation $\pi = 1342$ is identified with the matrix
$\left[\begin{smallmatrix} 1 & 0 & 0 & 0 \\
0 & 0 & 1 & 0 \\
0 & 0 & 0 & 1 \\
0 & 1 & 0 & 0
\end{smallmatrix}\right]$. Note this convention differs from that of most of the permutation patterns literature; it was chosen to be consistent with prior work on alternating sign matrices (see e.g.~\cite{ACG,JL}).

An \emph{alternating sign matrix} (ASM) is a square matrix where 
each entry is an element of $\{0,1,-1\}$,
every row and column sum to $1$, and
the non-zero entries of each row and column alternate in sign.
We let $\ASM_n$ denote the set of ASMs of size $n$, and we define $\ASM = \bigsqcup_{n\ge0} \ASM_n$. 
Note that, $S_n$ is a subset of $\ASM_n$. 

There are several conceivable ways of defining pattern containment for $\{0,\pm 1\}$ matrices, many of which are quite similar, so now we will carefully define the notion of pattern containment that we will be using in what follows. We call it \emph{classical containment} by analogy with the terminology in the permutation patterns literature. 

Let $B = {\left(B_{i,j}\right)}_{\substack{1 \le i \le m \\ 1 \le j \le n}}$ be an $m \times n$ matrix with entries from $\{0, \pm 1\}$, and let $A = {\left(A_{k,\ell}\right)}_{\substack{1 \le k \le p \\ 1 \le \ell \le q}}$ be a $p \times q$ matrix with entries from $\{0, \pm 1\}$. We say $B$ is \emph{classically contained in} $A$, or $A$ \emph{classically contains} $B$ (as a \emph{pattern}), if there are order-preserving injections $f \colon [m] \to [p]$ and $g \colon [n] \to [q]$ such that $B_{i,j} \not= 0$ implies $A_{f(i), g(j)} = B_{i,j}$. Equivalently, rows and columns can be deleted from $A$ to obtain an $m \times n$ matrix $A'$ such that every \emph{non-zero} entry of $B$ equals the corresponding entry of $A'$.

Under the assumption that $B$ is an $n \times n$ permutation matrix corresponding to $\pi \in S_n$, and that $A$ is a $p \times p$ ASM, which will be the case for most of this paper, the above definition can be restated as follows: a permutation $\pi$ is \emph{classically contained in} $A$ if there are order-preserving injections $f \colon [n] \to [p]$ and $g \colon [n] \to [p]$ such that $A_{f(i), g(\pi(i))} = 1$ for all $i$. Equivalently, rows and columns can be deleted from $A$ to obtain an $n \times n$ matrix with $1$'s in at least the positions where $B$ has $1$'s. This matches the definition in Johansson and Linusson \cite{JL}. An important observation is that, if $\chi(A)$ is the matrix obtained from $A$ by replacing every $-1$ with $0$, then the permutations classically contained in $A$ are the same as the permutations classically contained in $\chi(A)$; two $\{0, \pm1\}$ matrices with their $1$'s in the same positions will contain all the same permutations as patterns.

If $R$ is a set of permutations, then $S_n(R)$ (resp.\ $S(R)$) denotes the set of permutations in $S_n$ (resp.\ $S$) avoiding every element of $R$; similarly, $\ASM_n(R)$ (resp.\ $\ASM(R)$ denotes the set of ASMs in $\ASM_n$ (resp.\ $\ASM$) classically avoiding every element of $R$.

\section{Asymptotic dichotomy when avoiding one permutation pattern}
\label{sec:mainthm}


To prove Theorem \ref{thm:superexponential}, we start by introducing helpful notation. 
Let $I_k$ denote the $k \times k$ identity matrix (corresponding to the permutation $123 \ldots k$), and let $I_k'$ denote the matrix with $1$'s on the anti-diagonal and $0$'s everywhere else (corresponding to the permutation $k \ldots 321$). Let $0_k$ denote the $k \times k$ matrix of all $0$'s. 
Now let $n$ be a positive integer, and write $n=3k+r$ for $r\in \{0,1,2\}$. For each $\pi \in S_{k}$, we define two specific $n \times n$ matrices as follows. 
If $n$ is divisible by $3$, that is, if $r=0$, we define 
\[ \Gamma_n(\pi) = \begin{bmatrix}
0_{k} & I_{k} & 0_{k} \\
I_{k} & -\pi & I_{k} \\
0_{k} & I_{k} & 0_{k}
\end{bmatrix} \quad \text{and} \quad \Gamma_n'(\pi) = \begin{bmatrix}
0_{k} & I_{k}' & 0_{k} \\
I_{k}' & -\pi & I_{k}' \\
0_{k} & I_{k}' & 0_{k}
\end{bmatrix},
 \]
where $-\pi$ denotes the permutation matrix of $\pi$ with each $1$ replaced with $-1$. For example, writing only $+$ and $-$ for the entries $+1$ and $-1$, respectively, and leaving white spaces for the $0$ entries, we have 
\[
\Gamma_9(312) = \begin{bmatrix} {\scriptsize \begin{array}{ccc|ccc|ccc}
&&&+&&&&& \\
&&&&+&&&& \\
&&&&&+&&& \\
\hline
+&&&&&-&+&& \\
&+&&-&&&&+& \\
&&+&&-&&&&+ \\
\hline
&&&+&&&&& \\
&&&&+&&&& \\
&&&&&+&&&
\end{array} } \end{bmatrix}
\quad \text{and} \quad \Gamma_9'(312) = \begin{bmatrix} { \scriptsize \begin{array}{ccc|ccc|ccc}
&&&&&+&&& \\
&&&&+&&&& \\
&&&+&&&&& \\
\hline
&&+&&&-&&&+ \\
&+&&-&&&&+& \\
+&&&&-&&+&& \\
\hline
&&&&&+&&& \\
&&&&+&&&& \\
&&&+&&&&&
\end{array} } \end{bmatrix}. \]
If $r=1$ (resp. $r=2$) we modify the construction as follows, adding one (resp. two) row(s) and columns(s) having a $1$ entry at their intersection(s) and $0$ elsewhere: 

\[ \Gamma_n(\pi) = \begin{bmatrix}
0_{k} & I_{k} & 0_{k} & \mathbf{0}  \\
I_{k} & -\pi & I_{k} & \mathbf{0} \\
\mathbf{0} & \mathbf{0} & \mathbf{0} & 1 \\ 
0_{k} & I_{k} & 0_{k} & \mathbf{0}
\end{bmatrix} \quad \text{and} \quad \Gamma_n'(\pi) = \begin{bmatrix}
0_{k} & I_{k}' & 0_{k} & \mathbf{0}\\
\mathbf{0} &  \mathbf{0} & \mathbf{0} & 1 \\ 
I_{k}' & -\pi & I_{k}' & \mathbf{0} \\
0_{k} & I_{k}' & 0_{k}& \mathbf{0}
\end{bmatrix}
 \]

\[ \text{(resp. }  \Gamma_n(\pi) = \begin{bmatrix}
0_{k} & I_{k} & 0_{k} & \mathbf{0} & \mathbf{0}  \\
I_{k} & -\pi & I_{k} & \mathbf{0} & \mathbf{0} \\
\mathbf{0} &  \mathbf{0}& \mathbf{0} & 1 & 0  \\ 
\mathbf{0} & \mathbf{0} & \mathbf{0} & 0 & 1 \\ 
0_{k} & I_{k} & 0_{k} & \mathbf{0}& \mathbf{0} 
\end{bmatrix} \quad \text{and} \quad \Gamma_n'(\pi) = \begin{bmatrix}
0_{k} & I_{k}' & 0_{k} & \mathbf{0} & \mathbf{0}\\
\mathbf{0} & \mathbf{0} & \mathbf{0} & 0 & 1  \\ 
\mathbf{0} & \mathbf{0} & \mathbf{0} & 1 & 0 \\ 
I_{k}' & -\pi & I_{k}' & \mathbf{0}& \mathbf{0}\\
0_{k} & I_{k}' & 0_{k} & \mathbf{0} & \mathbf{0}
\end{bmatrix} \text{).}
 \]
 Each bold $\mathbf{0}$ in the above matrices indicates either a row or a column vector of $k$ $0$ entries, the direction of the vector being clear from the context. 
The matrices $\Gamma_n(\pi)$ and $\Gamma_n'(\pi)$ are ASMs; in fact, $\Gamma_n(\pi) \in \ASM(321)$ and $\Gamma_n'(\pi) \in \ASM(123)$.

Recall from the introduction that the set $X$ of \exceps~has been defined as 
\[X:=\{ 1, 12, 21, 132, 213, 231, 312, 2143, 2413, 3142, 3412\}.\]

\begin{MainTheorem} 
Suppose $\sigma$ is a permutation.  
If $\sigma\in X$, then there exists $c$ such that $|\ASM_n(\sigma)| \leq c^n$ for all $n$. 
Otherwise, $|\ASM_n(\sigma)| \ge \lfloor n/3 \rfloor !$ for all $n$.
\end{MainTheorem}

\begin{proof}
First observe that $\sigma \in X$ if and only if $\sigma$ avoids both $123$ and $321$.
Indeed, a permutation of size $5$ or larger cannot avoid both $123$ and $321$ (this is an instance of the Erd\H{o}s--Szekeres Theorem \cite{ES}) and going through all permutations of size at most $4$ then allows to check this claim. 

We start by proving the second statement ($|\ASM_n(\sigma)| \ge \lfloor n/3 \rfloor !$ for all $n$), in the case that $\sigma$ contains $321$. Under this assumption, for each $\pi \in S_{\lfloor n/3 \rfloor}$, it holds that $\sigma$ is not classically contained in $\Gamma_n(\pi)$, since $\Gamma_n(\pi)$ classically avoids $321$. Therefore, the set $\{\Gamma_n(\pi) \colon \pi \in S_{\lfloor n/3 \rfloor}\}$ comprises $\lfloor n/3 \rfloor!$ ASMs that classically avoid $\sigma$, proving that $|\ASM_n(\sigma)| \ge \lfloor n/3 \rfloor!$.

If $\sigma$ instead contains $123$, then the same proof works using $\Gamma_n'(\pi)$. 
With the characterization of $X$ given above, this completes the proof of the second statement. 

We turn to the proof of the first statement. 
An alternative description of the permutations in $X$ allows for a unified proof in this case. Indeed, we start by observing that the permutations in $X$ are exactly the ones whose elements can fit on a rectangle whose sides are exactly vertical and exactly horizontal. 
What we mean by that is that it is possible to identify an axis-parallel rectangle inside the permutation matrix such that all the $1$ entries lie on this rectangle. 
Proving that the permutations satisfying this property are exactly the permutations in $X$ is easily done by examining all permutations of size at most $4$ (size $4$ clearly being the maximum).
Specifically, permutations fitting on a rectangle must avoid both $123$ and $321$ and as such, they must be in $X$.

Let $\sigma \in X$ be an \excep. An ASM $A$ classically avoiding $\sigma$ is characterized by a pair of $\{0,1\}$ matrices $(A_+,A_-)$, where $A_+$ has been obtained from $A$ by replacing all $-1$ entries by $0$, and similarly $A_-$ has been obtained from $A$ by replacing all $1$ entries by $0$, and then replacing all $-1$ entries by $1$. 
Clearly, because $A$ classically avoids $\sigma$, so does $A_+$. 
We claim that it also holds that $A_-$ avoids $\sigma$. Indeed, because $\sigma$ is an \excep, its elements can be seen on an axis-parallel rectangle. As a consequence, if the $-1$ entries of $A$ would form a pattern $\sigma$, there must exist entries of $A$ equal to $1$  closer to the boundary of the matrix that also form a pattern $\sigma$, which contradicts that $A$ classically avoids $\sigma$. This is illustrated on Figure~\ref{F:excepIn1}. 

\begin{figure}[htbp]
    \begin{center}
\begin{tikzpicture}[scale=0.5]
\fill [blue!10!white] (0.5,1.5) rectangle (1.5,7.5);
\fill [blue!10!white] (1.5,6.5) rectangle (6.5,7.5);
\fill [blue!10!white] (5.5,1.5) rectangle (6.5,6.5);
\fill [red!10!white] (2.5,2.5) rectangle (4.5,3.5);
\fill [blue!10!white] (1.5,1.5) rectangle (5.5,2.5);
\fill [red!10!white] (1.5,4.5) rectangle (5.5,5.5);
\fill [red!10!white] (1.5,2.5) rectangle (2.5,4.5);
\fill [red!10!white] (4.5,2.5) rectangle (5.5,4.5);
\node at (1, 1) {$0$};
\node at (1, 2) {$0$};
\node at (1, 3) {\color{blue}{$1$}};
\node at (1, 4) {$0$};
\node at (1, 5) {$0$};
\node at (1, 6) {$0$};
\node at (1, 7) {$0$};
\node at (2, 1) {$0$};
\node at (2, 2) {$1$};
\node at (2, 3) {\color{red}{$-1$}};
\node at (2, 4) {$0$};
\node at (2, 5) {$1$};
\node at (2, 6) {$0$};
\node at (2, 7) {$0$};
\node at (3, 1) {$0$};
\node at (3, 2) {$0$};
\node at (3, 3) {$0$};
\node at (3, 4) {$1$};
\node at (3, 5) {$0$};
\node at (3, 6) {$0$};
\node at (3, 7) {$0$};
\node at (4, 1) {$0$};
\node at (4, 2) {$0$};
\node at (4, 3) {$1$};
\node at (4, 4) {$0$};
\node at (4, 5) {\color{red}{$-1$}};
\node at (4, 6) {$0$};
\node at (4, 7) {\color{blue}{$1$}};
\node at (5, 1) {$1$};
\node at (5, 2) {$0$};
\node at (5, 3) {$0$};
\node at (5, 4) {\color{red}{$-1$}};
\node at (5, 5) {$0$};
\node at (5, 6) {$1$};
\node at (5, 7) {$0$};
\node at (6, 1) {$0$};
\node at (6, 2) {$0$};
\node at (6, 3) {$0$};
\node at (6, 4) {\color{blue}{$1$}};
\node at (6, 5) {$0$};
\node at (6, 6) {$0$};
\node at (6, 7) {$0$};
\node at (7, 1) {$0$};
\node at (7, 2) {$0$};
\node at (7, 3) {$0$};
\node at (7, 4) {$0$};
\node at (7, 5) {$1$};
\node at (7, 6) {$0$};
\node at (7, 7) {$0$};
\end{tikzpicture}
    \end{center}
    \caption{An \excep~among the $-1$'s of an ASM (shown in red) forces an \excep~among its $1$'s (shown in blue). }
    \label{F:excepIn1}
\end{figure}

Therefore, letting $M_n(\sigma)$ denote the set of $\{0,1\}$ matrices avoiding $\sigma$, we have $|\ASM_n(\sigma)| < |M_n(\sigma)|^2$.
By the Marcus--Tardos theorem~\cite[Thm.\ 9]{MarcusTardos} we have $|M_n(\sigma)| \le 15^{2k^4 \binom{k^2}{k}n}$, where $k$ is the size of $\sigma$.
It follows that $|\ASM_n(\sigma)|$ is at most $15^{4k^4 \binom{k^2}{k}n}$, concluding the proof of our statement.
\end{proof}

\begin{remark}

It would also be interesting to enumerate $\ASM(321)$, asymptotically if not exactly. Indeed, every ASM class avoiding one pattern not in $X$ contains either $\ASM(321)$ or $\ASM(123)$, which are symmetries of each other. Thus, a lower bound on $\ASM(321)$ would also be a lower bound to replace $\lfloor n/3 \rfloor !$ in Theorem \ref{thm:superexponential}.
\end{remark}

\begin{remark}
In the case where $\sigma$ contains $321$, we can construct more elements of $\ASM_n(\sigma)$ using the $\Gamma(\pi)$ ASMs as blocks: for any choice of positive integers $n_1, \dots, n_l$ such that $n_1 + \dots + n_l = n$, and any choice of ASMs of the form $\Gamma_{n_1}(\pi^{(1)}), \Gamma_{n_2}(\pi^{(2)}), \dots, \Gamma_{n_l}(\pi^{(l)})$ (where $\pi^{(1)}, \dots, \pi^{(l)}$ are permutations of the appropriate sizes), we can form the block-diagonal matrix
\[ \begin{bmatrix} \Gamma_{n_1}(\pi^{(1)}) & & & \\
& \Gamma_{n_2}(\pi^{(2)}) & & \\
& & \ddots & \\
& & & \Gamma_{n_l}(\pi^{(l)})
\end{bmatrix}, \]
which is in $\ASM_n(321)$ and hence in $\ASM_n(\sigma)$. The number of ASMs that can be formed in this way is
\[ \sum_{\substack{n_1,\dots,n_l >0 \\ n_1+\dots+n_l = n}} \lfloor n_1/3 \rfloor !\, \lfloor n_2/3 \rfloor! \, \dots \, \lfloor n_l/3 \rfloor! . \]
The resulting numbers are asymptotically a constant multiple of $\lfloor n/3 \rfloor !$ (this constant depends on $n$ mod $3$), an improvement which to us is not significant enough to be worth complicating our theorem.
\end{remark}

We now turn to the question of enumerating $\ASM(R)$ for a set of permutations $R$. If one of the eleven patterns in $X$ is an element of $R$, then the counting sequence is bounded above by an exponential, 
by Theorem~\ref{thm:superexponential}. On the other hand, if every element of $R$ contains $123$, or if every element of $R$ contains $321$, then the proof of Theorem \ref{thm:superexponential} shows that $|\ASM_n(R)| \ge \lfloor n/3 \rfloor!$ in this case. We can make this last statement even stronger, with Theorem \ref{thm:incdec} below. 

Let $\mathcal{C}^\text{inc}$ (resp.\ $\mathcal{C}^\text{dec}$) be the class of permutations $\pi$ such that the matrix of $\pi$ has a block form $\begin{pmatrix} 0 & A & 0 \\ B & 0 & C \\ 0 & D & 0 \end{pmatrix}$ in which the entries in each block $A$, $B$, $C$, and $D$ form an increasing (resp.\ decreasing) pattern. 
For those readers familiar with (monotone) grid classes (defined in~\cite{grid}), we note that $\mathcal{C}^\text{inc} = \text{Grid} \left( \begin{smallmatrix}
0 & 1 & 0 \\
1 & 0 & 1 \\
0 & 1 & 0
\end{smallmatrix} \right)$ and $\mathcal{C}^\text{dec} = \text{Grid} \left( \begin{smallmatrix}
0 & -1 & 0 \\
-1 & 0 & -1 \\
0 & -1 & 0
\end{smallmatrix} \right)$. 
The following theorem has the same proof as Theorem \ref{thm:superexponential}:

\begin{thm} \label{thm:incdec}
Let $R$ be a set of permutations. If $R \subseteq S \smallsetminus \mathcal{C}^\textnormal{inc}$ or $R \subseteq S \smallsetminus \mathcal{C}^\textnormal{dec}$,
then $|\ASM_n(R)| \ge \lfloor (n/3) \rfloor!$.
\end{thm}

\begin{proof}
First assume $R \subseteq S \smallsetminus \mathcal{C}^\text{inc}$. For each $n$, the set $\{\Gamma_n(\pi) \colon \pi \in S_{\lfloor n/3 \rfloor}\}$ comprises $\lfloor n/3 \rfloor!$ ASMs such that every permutation pattern contained in these $\lfloor n/3 \rfloor!$ ASMs is an element of $\mathcal{C}^\text{inc}$. Hence these $\lfloor n/3 \rfloor!$ ASMs are in $\ASM(R)$. Thus $|\ASM_n(R)| \ge \lfloor n/3 \rfloor!$.

If $R \subseteq S \smallsetminus \mathcal{C}^\text{dec}$, then the same proof works using $\Gamma'_n(\pi)$.
\end{proof}

Note that $\mathcal{C}^\text{inc}$ and $\mathcal{C}^\text{dec}$ are contained in $S(321)$ and $S(123)$ respectively, because every permutation in $\mathcal{C}^\text{inc}$ (resp.\ $\mathcal{C}^\text{dec}$) is a union of two increasing (resp.\ decreasing) sequences, and such a permutation cannot contain a $321$ (resp.\ $123$) pattern. Thus Theorem \ref{thm:incdec} actually implies the super-exponential part of Theorem \ref{thm:superexponential}.

Finally, we show below that for any permutation $\sigma$, and for any fixed integer $k$, the sequence enumerating $n\times n$ ASMs that classically avoid $\sigma$ and contain exactly $k$ negative ones grows at most exponentially with $n$. 
This is in interesting contrast with the super-exponential lower bound on $|\ASM_n(\sigma)|$ when $\sigma \notin X$. Recall  $\ASM_{n,k}(\sigma)$ denotes the set of ASMs in $\ASM_n(\sigma)$ containing exactly $k$ negative ones.

\begin{kNegativeOnesThm}
    Suppose $\sigma$ is a permutation and $k$ is a non-negative integer.  
Then there exists $c$ such that $|\ASM_{n,k}(\sigma)| \leq c^n$ for all $n$. 
\end{kNegativeOnesThm}

\begin{proof}
The proof relies on the same ideas used in the proof of the exponential upper bound part of Theorem~\ref{thm:superexponential}. Indeed, the number of ways to put ones in such an ASM is bounded above by an exponential, again by the Marcus--Tardos theorem.  And the number of ways to put negative ones in such an ASM is bounded above by $\binom{n^2}{k}$, which is a polynomial of degree $2k$ in $n$. So the number of ASMs we are counting is bounded above by the product of an exponential and a polynomial, which is bounded above by an exponential.
\end{proof}

\section{An explicit upper bound for $|\ASM_n(3412)|$}
\label{sec:abetterbound}
Theorem \ref{thm:superexponential} specifies all of the patterns $\sigma$ such that $\ASM(\sigma)$ is at most exponential. Since there are only $11$ patterns, it would be nice to exactly enumerate or give an explicit bound for each pattern class (and a better one than the one derived through the Marcus--Tardos theorem and our proof of Theorem \ref{thm:superexponential}). The enumeration of $\ASM(\sigma)$ is trivial for $\sigma \in \{1, 12, 21\}$; the counting sequence for each of these is eventually $0$ or eventually $1$. The exact enumeration of $\ASM(132)$, which by symmetry is equal to the enumeration of $\ASM(\sigma)$ for $\sigma \in \{213, 231, 312\}$, was completed by Johansson and Linusson \cite{JL}. 
This leaves two other cases (up to symmetry; specifically, taking the reverse): 
\[ |\ASM_n(3412)|=|\ASM_n(2143)|  \quad \text{and} \quad |\ASM_n(3142)|=|\ASM_n(2413)|. \]
We study the case of 3412 (or, by symmetry, 2143) in Theorem~\ref{T:3412} below, giving a better bound than the one that can be derived from 
Theorem \ref{thm:superexponential} and its proof.

We begin with the following lemma.
\begin{lem} \label{lem:2143}
Let $A$ be an ASM. If $A$ has a $-1$ strictly southwest of another $-1$, then $A$ classically contains the pattern $3412$. If $A$ has a $-1$ strictly northwest of another $-1$, then $A$  classically contains the pattern $2143$. 
\end{lem}
\begin{proof}
Let $A\in \ASM_n$ and $(i,j)$ and $(i',j')$ be the positions of two $(-1)$'s in $A$, and suppose that $(i,j)$ is strictly southwest of $(i',j')$ --- that is, $i > i'$ and $j < j'$. Then there must be a $1$ strictly to the north of $(i',j')$ and another $1$ strictly to the east of $(i',j')$ forming the $3,4$ respectively of the $3412$ pattern.  Then the $1$ strictly to the west of $(i,j)$ and the $1$ strictly to the south of $(i,j)$ form the $1,2$ of the $3412$ pattern.  See Figure~\ref{F:TwoNegOne1} for an example of this.         

The proof for $2143$ is a symmetry of the above.
\end{proof}

\begin{figure}[htbp]
    \begin{center}

        $\begin{pmatrix}
            0 & 0 & 0 & \boxed{1} & 0 & 0 \\
            0 & 0 & 0 & 0 & 0 & 1 \\
            0 & 1 & 0 & -1 & \boxed{1} & 0\\
            \boxed{1} & -1 & 1 & 0 & 0 & 0 \\
            0 & 0 & 0 & 1 & 0 & 0 \\
            0 & \boxed{1} & 0 & 0 & 0 & 0 \\
        \end{pmatrix}   $
    
    \end{center}
    \caption{An illustration of Lemma 4.1:  An ASM with two $(-1)$s where one is strictly southwest of the other will have a 3412 pattern.}
    \label{F:TwoNegOne1}
\end{figure}

\begin{BetterBoundThm}
The number of alternating sign matrices that classically avoid $3412$ has the following exponential upper bound:
$|\ASM_n(3412)| \leq 62208^n$. The same holds for  $\ASM_n(2143)$.
\end{BetterBoundThm}
\begin{proof}
Since an ASM classically avoids $2143$ if and only if its reverse classically avoids $3412$, we have $|\ASM_n(2143)| = |\ASM_n(3412)|$. We proceed to prove the theorem for $\ASM_n(3412)$.

Let $A \in \ASM_n(3412)$ and let $k$ be the number of negative ones in $A$. Let $(i_1, j_1), \ldots, (i_k, j_k)$ be the positions of the $(-1)$'s in $A$. Lemma~\ref{lem:2143} indicates that these positions appear in a linear order from the northwest to the southeast. More precisely, Lemma~\ref{lem:2143} implies that the positions $(i_1, j_1), \ldots, (i_k, j_k)$ can be totally ordered such that if $a \le b$ then $i_a \le i_b$ and $j_a \le j_b$. Because two $(-1)$'s in the same row or column need to have a $1$ between them, this implies in particular that $k\leq n$. 

In addition, the positions of the $(-1)$'s are therefore determined by choosing subsets $R \subseteq [n]$ and $C \subseteq [n]$ describing the rows and columns containing a $(-1)$, together with a sequence of $S$ (south), $E$ (east), and $D$ (diagonal) steps, of $k-1$ steps in total, indicating the direction to follow to find the next $(-1)$ in the total order defined above. The number of such sequences of $S$, $E$ and $D$ is at most $3^k \leq 3^n$. Therefore, the number of possible placements of the $(-1)$'s is at most
\[ 2^n \cdot 2^n \cdot 3^n = 12^n. \]
Now we count the ways of placing the $1$'s for a given placement of $(-1)$'s. First, we create a path of adjacent matrix entries as follows. This path contains all the positions $(i_1, j_1), \ldots, (i_k, j_k)$ with $(-1)$'s (defined above), as well as all matrix entries between two $(-1)$'s in the same row or between two $(-1)$'s in the same column. We also include $(i_1,j')$ for all $j'<j_1$ and $(i_k,j'')$ for all $j''>j_k$. Finally, for each pair $(i_a,j_a),(i_{a+1},j_{a+1})$ of $(-1)$'s that form a $D$ in the sequence (i.e.\ $(i_{a+1}, j_{a+1})$ is southeast of $(i_a, j_a)$), include all $(i_a,j')$ with $j_a<j'\leq j_{a+1}$ and all $(i',j_{a+1})$ with $i_a\leq i' < i_{a+1}$. See Figure~\ref{F:Path} for an illustration. 

\begin{figure}[htbp]
    \begin{center}

        $\begin{pmatrix}
            
            0 & 0 & 1 & 0 & 0 & 0 & 0 & 0 & 0 & 0 & 0 & 0 & 0 & 0 & 0 \\
            \boxed{1} & \boxed{0} & \boxed{-1} & \boxed{0} & \boxed{0} & 0 & 0 & 1 & 0 & 0 & 0 & 0 & 0 & 0 & 0 \\
            0 & 0 & 0 & 0 & \boxed{1} & 0 & 0 & 0 & 0 & 0 & 0 & 0 & 0 & 0 & 0 \\
            0 & 0 & 1 & 0 & \boxed{0} & 0 & 0 & 0 & 0 & 0 & 0 & 0 & 0 & 0 & 0 \\
            0 & 1 & 0 & 0 & \boxed{-1} & \boxed{0} & \boxed{1} & \boxed{-1} & \boxed{0} & \boxed{1} & 0 & 0 & 0 & 0 & 0 \\
            0 & 0 & 0 & 0 & 0 & 0 & 0 & 0 & 0 & \boxed{0} & 0 & 1 & 0 & 0 & 0 \\
            0 & 0 & 0 & 0 & 0 & 0 & 0 & 0 & 1 & \boxed{-1} & 1 & 0 & 0 & 0 & 0 \\
            0 & 0 & 0 & 0 & 0 & 1 & 0 & 0 & 0 & \boxed{0} & 0 & 0 & 0 & 0 & 0 \\
            0 & 0 & 0 & 0 & 0 & 0 & 0 & 0 & 0 & \boxed{1} & 0 & 0 & 0 & 0 & 0 \\
            0 & 0 & 0 & 0 & 0 & 0 & 0 & 0 & 0 & \boxed{0} & 0 & 0 & 1 & 0 & 0 \\
            0 & 0 & 0 & 0 & 0 & 0 & 0 & 1 & 0 & \boxed{-1} & \boxed{0} & \boxed{0} & \boxed{0} & 0 & 1 \\
            0 & 0 & 0 & 0 & 0 & 0 & 0 & 0 & 0 & 1 & 0 & 0 & \boxed{0} & 0 & 0 \\
            0 & 0 & 0 & 0 & 1 & 0 & 0 & 0 & 0 & 0 & 0 & 0 & \boxed{-1} & \boxed{1} & \boxed{0} \\
            0 & 0 & 0 & 0 & 0 & 0 & 0 & 0 & 0 & 0 & 0 & 0 & 1 & 0 & 0 \\
            0 & 0 & 0 & 1 & 0 & 0 & 0 & 0 & 0 & 0 & 0 & 0 & 0 & 0 & 0 
        \end{pmatrix}   $
    
    \end{center}
    \caption{An illustration of the path of adjacent matrix entries defined in Theorem~\ref{T:3412}.}
    \label{F:Path}
\end{figure}

Each $1$ is either on this path, strictly above the path, or strictly below the path. The number of ways of placing the $1$'s that are on the path is at most $2^{2n}$, because the path has at most $2n$ entries. The $1$'s strictly above the path form a $3412$-avoiding permutation on some subset of the rows and columns, so the number of ways of placing them is at most $2^n \cdot 2^n \cdot |S_n(3412)|$; and the same enumeration holds for the $1$'s strictly below the path. So the total number of ways of placing the $1$'s is
\[ 2^{2n} \cdot [2^n \cdot 2^n \cdot |S_n(3412)|]^2 = 64^n \cdot |S_n(3412)|^2 \le 64^n \cdot (9^n)^2 = 5184^n \]
for $n$ large enough, since $|S_n(3412)| = |S_n(4321)| \le 9^n$ \cite[Thm.\ 4.10 \& Thm.\ 4.12]{Bona}. Therefore, $|\ASM_n(3412)| \le 12^n \cdot 5184^n = 62208^n$.
\end{proof}

Theorem~\ref{T:3412} provides an exponential upper bound for $|\ASM_n(3412)|$, but it does not shed any light on the actual asymptotic formula of this sequence. 
We believe it would be interesting to further study ASMs classically avoiding $3412$, to possibly determine their exact or asymptotic enumeration. Based on the actual numbers for small $n$, our upper bound of $62208^n$ appears far from optimal, but we note that it is still much better than the upper bound from Theorem~\ref{thm:superexponential} and its proof, which is exponential with base $15^{1863680}$.

We end this section with a lemma that could possibly be useful in this direction, as it identifies additional structure in $3412$-avoiding ASMs.

\begin{lem}~\label{L:stay_on _path}
    If an ASM that classically avoids $3412$ has two $(-1)$s where one $-1$ in position $(i_a,j_a)$ is strictly northwest of the other $-1$ in position $(i_b,j_b)$, that is $i_a < i_b$ and $j_a < j_b$, then there must be a $1$ on at least one of the two paths $P_1$ and $P_2$, where path $P_1$ goes from $(i_a,j_a)$ to $(i_b,j_b)$ by a sequence of $E = (0,1)$ steps followed by a sequence of $S = (1,0)$ steps and path $P_2$ goes from $(i_a,j_a)$ to $(i_b,j_b)$ by a sequence of $S$ steps followed by a sequence of $E$ steps.
\end{lem}

\begin{proof}
    Suppose there are no $1$s along either of $P_1,P_2$.  Then there are $1$s in the following places:
    \begin{itemize}
        \item north of $(i_b,j_b)$, which since it is not on $P_1$, is also strictly northeast of $(i_a,j_a)$ in position $(i_a -y_1,j_b)$ for some $y_1  \geq 1$.
        \item east of $(i_a,j_a)$,  which since it is not on $P_1$, is also strictly northeast of $(i_b,j_b)$ in position $(i_a,j_b + x_1)$ for some $x_1  \geq 1$.
        \item west of $(i_b,j_b)$, which since it is not on $P_2$, is also strictly southwest of  $(i_a,j_a)$ in position $(i_b,j_a -x_2)$ for some $x_2  \geq 1$.
        \item south of $(i_a,j_a)$, which since it is not on $P_2$, is also strictly southwest of  $(i_b,j_b)$ in position $(i_b + y_2,j_a)$ for some $y_2  \geq 1$.
    \end{itemize}
    Since $i_a -y_1 < i_a < i_b  < i_b +y_2$ and $j_a-x_2 < j_a < j_b < j_b + x_1$, the $1$s in the above listed coordinates form a $3412$ pattern.  See Figure~\ref{F:stay_on_path} for an illustration.
\end{proof}

\begin{figure}[htbp]
    \begin{center}

        $\begin{pmatrix}
            0 & 0 & 0 & 0 & \boxed{1} & 0 & 0 \\
            0 & 0 & 1 & 0 & 0 & 0 & 0 \\
            1 & 0 &-1 & 0 & 0 & 0 &\boxed{1} \\
            0 & 0 & 0 & 1 & 0 & 0 & 0 \\
            0 & \boxed{1} & 0 & 0 & -1 & 1 & 0 \\
            0 & 0 & \boxed{1} & 0 & 0 & 0 & 0 \\
            0 & 0 & 0 & 0 & 1 & 0 & 0
        \end{pmatrix}   $
    
    \end{center}
    \caption{An illustration of the proof of Lemma~\ref{L:stay_on _path}.}
    \label{F:stay_on_path}
\end{figure}

\section{Toward an exact enumeration of  $\ASM_n(2143,3412)$}
\label{sec:toward}
In this section, we prove some lemmas which may be useful toward finding an exact enumeration of $\ASM_n(2143,3412)$. Using Sage~\cite{Sage}, we have  obtained the first few terms of the sequence for $|\ASM_n(2143,3412)|$: $1, 2, 7, 38, 228, 1232, 5888$. In particular, we show that any such ASM has at most three negative ones. Our main result of this section is  Theorem~\ref{thm:exactly3enumeration} enumerating the number of ASMs in $\ASM_n(2143,3412)$ with exactly three negative ones. 

\begin{lem} \label{lem:allinarow}
For any $A\in \ASM_n(2143,3412)$, all the negative ones are in a single row or a single column. 
\end{lem}
\begin{proof} This follows directly from Lemma~\ref{lem:2143}, as an ASM containing two $-1$s that are not in the same row or column will imply the ASM has a $2143$ pattern or a $3412$ pattern.
\end{proof}

\begin{lem}
\label{lem:atmost3}
Any $A\in \ASM_n(2143,3412)$ has at most $3$ negative ones.
\end{lem}
\begin{proof}
Suppose $A\in \ASM_n(2143,3412)$ has 4 or more negative ones. By Lemma~\ref{lem:allinarow}, these negative ones lie in the same row or column. Suppose  they are in a row. So the row begins with the following sign pattern, omitting the zeros: $+ - + - + - + - +$. Now there must be $+1$'s in the same column above and below each $-1$. 

Consider the $1$ due north of the first $-1$ in the row and the $1$ at the beginning of the row whose column values form a $21$ pattern.  Then the $1$s to the south of the remaining  $-1$s in the row must form an increasing sequence of column values as the row positions increase to avoid creating a $2143$ pattern.  Specifically, the $1$s to the south of the second and third  $-1$s in the row must form an increasing sequence of column values as shown in Figure~\ref{F:fourNegOnes} (left). 
\begin{figure}[htbp]
    \begin{center}
        $\begin{pmatrix}
            0 & \boxed{1} & 0 & 0 & 0 & 0 & 0 & 0 & 0\\
            \vdots & \vdots & \vdots & \vdots & \vdots & \vdots & \vdots & \vdots & \vdots \\
            \boxed{1} & -1 & 1 & -1 & 1 & -1 & 1 & -1 & 1\\
            \vdots & \vdots & \vdots & \vdots & \vdots & \vdots & \vdots & \vdots & \vdots \\
            0 & 0 & 0 & \boxed{1} & 0 & 0 & 0 & 0 & 0\\
            \vdots & \vdots & \vdots & \vdots & \vdots & \vdots & \vdots & \vdots & \vdots \\
            0 & 0 & 0 & 0 & 0 & \boxed{1} & 0 & 0 & 0\\
        \end{pmatrix}   $ \qquad \qquad $\begin{pmatrix}
            0 & 0 & 0 & 0 & 0 & 0 & 0 & \boxed{1} & 0\\
            \vdots & \vdots & \vdots & \vdots & \vdots & \vdots & \vdots & \vdots & \vdots \\
            1 & -1 & 1 & -1 & 1 & -1 & 1 & -1 & \boxed{1}\\
            \vdots & \vdots & \vdots & \vdots & \vdots & \vdots & \vdots & \vdots & \vdots \\
            0 & 0 & 0 & \boxed{1} & 0 & 0 & 0 & 0 & 0\\
            \vdots & \vdots & \vdots & \vdots & \vdots & \vdots & \vdots & \vdots & \vdots \\
            0 & 0 & 0 & 0 & 0 & \boxed{1} & 0 & 0 & 0\\
        \end{pmatrix}   $
    \end{center}
    \caption{The configuration of a submatrix in an ASM with four $-1$s in the same row classically avoiding $2143$ (left) necessarily contains $3412$ (right).}
    \label{F:fourNegOnes}
\end{figure}

However, then the $1$ due north of the last $-1$ in the row and the $1$ at the end of the row  form a $3412$ pattern with the $1$s to the south of the second and third  $-1$s, as shown in Figure~\ref{F:fourNegOnes} (right). 
\end{proof}

Since by the above lemma, there are at most three negative ones, we may examine each case. The case with no negative ones corresponds to permutations avoiding $2143$ and $3412$. These are also known as \emph{skew-merged} permutations, first studied by Stankova~\cite{Stankova} and enumerated by Atkinson~\cite{Atkinson}. 
The name skew-merged permutation echoes a geometrical description of these permutations, proved to be equivalent to the avoidance of $2143$ and $3412$: 
namely, skew-merged permutations can be described as merges of one increasing sequence with a decreasing sequence. In the permutation matrix, this means that the $1$ entries form an $X$ shape. For example, the permutation $ 2 \, 5\,  10\,  8 \, 6\,  7\,  9\,  11\,  4\,  3\,  12 \, 1 $ avoids $2143$ and $3412$, and is a merge of the increasing sequence $ 2 \, 5\,  6\,  7\,  9\,  11\,  12$  with the decreasing sequence $10\,  8 \, 4\,  3\, 1 $. Its permutation matrix shown in Figure~\ref{F:skewmerge} exemplifies the $X$ shape, with the increasing sequence in red and the decreasing sequence in blue. 

\begin{figure}[htbp]
    \begin{center}
        $\begin{pmatrix}
            0 & {\color{red}{\mathbf{1}}} & 0 & 0 & 0 & 0 & 0 & 0 & 0 & 0 & 0 & 0\\
            0 & 0 & 0 & 0 & {\color{red}{\mathbf{1}}} & 0 & 0 & 0 & 0 & 0 & 0 & 0\\
            0 & 0 & 0 & 0 & 0 & 0 & 0 & 0 & 0 & {\color{blue}{\mathbf{1}}} & 0 & 0\\
            0 & 0 & 0 & 0 & 0 & 0 & 0 & {\color{blue}{\mathbf{1}}} & 0 & 0 & 0 & 0\\
            0 & 0 & 0 & 0 & 0 & {\color{red}{\mathbf{1}}} & 0 & 0 & 0 & 0 & 0 & 0\\
            0 & 0 & 0 & 0 & 0 & 0 & {\color{red}{\mathbf{1}}} & 0 & 0 & 0 & 0 & 0\\
            0 & 0 & 0 & 0 & 0 & 0 & 0 & 0 & {\color{red}{\mathbf{1}}} & 0 & 0 & 0\\
            0 & 0 & 0 & 0 & 0 & 0 & 0 & 0 & 0 & 0 & {\color{red}{\mathbf{1}}} & 0\\
            0 & 0 & 0 & {\color{blue}{\mathbf{1}}} & 0 & 0 & 0 & 0 & 0 & 0 & 0 & 0\\
            0 & 0 & {\color{blue}{\mathbf{1}}} & 0 & 0 & 0 & 0 & 0 & 0 & 0 & 0 & 0\\
            0 & 0 & 0 & 0 & 0 & 0 & 0 & 0 & 0 & 0 & 0 & {\color{red}{\mathbf{1}}}\\
            {\color{blue}{\mathbf{1}}} & 0 & 0 & 0 & 0 & 0 & 0 & 0 & 0 & 0 & 0 & 0\\
        \end{pmatrix}   $
    \end{center}
    \caption{Example of a skew-merged permutation.}
    \label{F:skewmerge}
\end{figure}

\begin{thm}[\cite{Atkinson}]  
The number of ASMs in  $\ASM_n(2143,3412)$ with zero negative ones, that is the number of skew-merged permutation matrices of size $n$, is
$$s_n = \binom{2n}{n} - \sum_{m=0}^{n-1} 2^{n-m-1}\binom{2m}{m}$$ 
which has generating function
$$\sum_{n=0}^{\infty}s_nz^n = \frac{(1-3z)}{(1-2z)\sqrt{1-4z}}.$$
\end{thm}

For the remainder of this section, we focus on the set of ASMs in $\ASM_n(2143,3412)$ with exactly three negative ones. 
Because all of the rotations in the dihedral group $D_4$ of the pattern class $\{2143,3412\}$ return the same class, in our arguments, we will assume without loss of generality that the three $-1$s 
are in the same row $i$. Specifically, say that the nonzero entries in row $i$ are in columns $a,b,c,d,e,f,g$ where $a<b<c<d<e<f<g$. In particular, there is a $1$ (resp. $-1$) entry in each column $a,c,e$ and $g$ (resp. $b,d$ and $f$) of row~$i$.

\begin{lem} 
\label{lem:iplusminus1}
Let $A \in \ASM_n(2143,3412)$ with exactly three $-1$s in row $i$. Then $A$ will have the $1$s above and below the center $-1$ (i.e.\ in column $d$) in rows $i-1$ and $i+1$.
\end{lem}
\begin{proof}
    First we show the $1$ that appears above the center $-1$ in row $i$ must in fact be in row $i-1$.  Otherwise, there is a $1$ in row $i-1$ that is either
    \begin{enumerate}
        \item in an earlier column than the $1$ above the center $-1$, so that these two $1$s form a $21$ pattern with the last $1$ of row $i$ and the $1$ below the eastmost $-1$, completing the forbidden $2143$ pattern as show in Figure~\ref{F:3neg1s_closemiddle} (left) or
        \item in a later column than the $1$ above the center $-1$, so that these two $1$s form a $34$ pattern with first $1$ of row $i$ and the $1$ below the westmost $-1$, completing the forbidden $3412$ pattern as show in Figure~\ref{F:3neg1s_closemiddle} (right).
    \end{enumerate}
    \begin{figure}[htbp]
    \begin{center}
        $\begin{pmatrix}
            0 & 0 & 0 & \boxed{1} & 0 & 0 & 0 \\
            \vdots & \vdots & \vdots & \vdots & \vdots & \vdots &  \vdots \\
            0 & \boxed{1} & 0 & 0 & 0 & 0 & 0 \\
            \vdots & \vdots & \vdots & \vdots & \vdots & \vdots &  \vdots \\
            1 & -1 & 1 & -1 & 1 & -1 & \boxed{1} & \\
            \vdots & \vdots & \vdots & \vdots & \vdots & \vdots &  \vdots \\
            0 & 0 & 0 & 0 & 0 & \boxed{1} & 0 \\
        \end{pmatrix}   $ \qquad \qquad         $\begin{pmatrix}
            0 & 0 & 0 & \boxed{1} & 0 & 0 & 0 \\
            \vdots & \vdots & \vdots & \vdots & \vdots & \vdots &  \vdots \\
            0 & 0 & 0 & 0 & 0 & \boxed{1} & 0 \\
            \vdots & \vdots & \vdots & \vdots & \vdots & \vdots &  \vdots \\
            \boxed{1} & -1 & 1 & -1 & 1 & -1 & 1 & \\
            \vdots & \vdots & \vdots & \vdots & \vdots & \vdots &  \vdots \\
            0 & \boxed{1} & 0 & 0 & 0 & 0 & 0 \\
        \end{pmatrix}   $
    \end{center}
    \caption{Part of an ASM with three $-1$s in the same row where there is a $1$ to the southwest or southeast of the $1$ above the center $-1$.}
    \label{F:3neg1s_closemiddle}
    \end{figure}

    Similarly, the $1$ below the center $-1$ must appear exactly in row $i+1$.  
    \end{proof}

    \begin{lem}
    \label{lem:abcdefg_adjacent}
        Let $A \in \ASM_n(2143,3412)$ with exactly three $-1$s in row $i$. Then the non-zero entries in row $i$ all appear in adjacent columns.
    \end{lem}
    \begin{proof}
     Suppose by way of contradiction that the row $i$ of $A$ has at least one $0$ between two nonzero entries.  If so, then there must be at least four $1$s above row $i$ or at least four $1$s below row $i$ in the columns strictly between $a$ and $g$.
     
     Suppose there are four $1$s above row $i$ in the columns strictly between $a$ and $g$ 
     (the other case is symmetric).  Then there will be two $1$s in the columns to the west of column $d$ (that of the center $-1$) or two $1$s to the east of column $d$.  Again, the symmetry makes these cases similar, so assume that there are two $1$s in the columns to the west of column $d$: the one above the $-1$ in column $b$, and an additional one, say in column $\beta$. Then we have three options:
     \begin{enumerate}
         \item The two $1$s in columns $b$ and $\beta$ form a $21$ pattern.  Then the eastmost $1$ of row $i$ (in column $g$) and the $1$ below the center $-1$ (in column $d$) complete a $2143$ pattern as shown in Figure~\ref{F:3neg1s_gapsA}.
    \begin{figure}[htbp]
    \begin{center}
        $\begin{pmatrix}
            0 & 0 & \boxed{1} & 0 & 0 & 0 & 0 & 0\\
            \vdots & \vdots & \vdots & \vdots & \vdots & \vdots &  \vdots &  \vdots\\
            0 & \boxed{1} & 0 & 0 & 0 & 0 & 0 & 0\\
            \vdots & \vdots & \vdots & \vdots & \vdots & \vdots &  \vdots &  \vdots \\
            1 & -1 & 0 & 1 & -1 & 1 & -1 & \boxed{1} & \\
            0 & 0 & 0 & 0 & \boxed{1} & 0 & 0 & 0\\
        \end{pmatrix}   $
    \end{center}
    \caption{Part of an ASM with three $-1$s in the same row where there is a $0$ between two nonzero entries in row $i$ left of center  and the $1$ in that column is above row $i$ forming a $21$ pattern with the $1$ above the westmost $-1$.}
    \label{F:3neg1s_gapsA}
    \end{figure}
    \item The two $1$s in columns $b$ and $\beta$ form a $12$ pattern and $\beta < b$.  Then the $1$ above the $0$ in row $i$ (in column $\beta$), the westmost $1$ of row $i$ (in column $a$), 
    the $1$ below the center $-1$ (in column $d$), and the $1$ below the westmost $-1$ (in column $b$) 
    form a $2143$ pattern as shown in Figure~\ref{F:3neg1s_gaps} (left).
    \begin{figure}[htbp]
    \begin{center}
        $\begin{pmatrix}
            0 & \boxed{1} & 0 & 0 & 0 & 0 & 0 & 0\\
            \vdots & \vdots & \vdots & \vdots & \vdots & \vdots &  \vdots &  \vdots\\
            0 & 0 & 1 & 0 & 0 & 0 & 0 & 0\\
            \vdots & \vdots & \vdots & \vdots & \vdots & \vdots &  \vdots &  \vdots \\
            \boxed{1} & 0 & -1 & 1 & -1 & 1 & -1 & 1 & \\
             0 & 0 & 0 & 0 & \boxed{1} & 0 & 0 & 1 & \\
            \vdots & \vdots & \vdots & \vdots & \vdots & \vdots &  \vdots &  \vdots\\
            0 & 0 & \boxed{1} & 0 & 0 & 0 & 0 & 0\\
        \end{pmatrix}   $ \qquad \qquad  $\begin{pmatrix}
            0 & 1 & 0 & 0 & 0 & 0 & 0 & 0 \\
            \vdots & \vdots & \vdots & \vdots & \vdots & \vdots &  \vdots &  \vdots\\
            0 & 0 & \boxed{1} & 0 & 0 & 0 & 0 & 0\\
            \vdots & \vdots & \vdots & \vdots & \vdots & \vdots &  \vdots &  \vdots \\
            0 & 0 & 0 & 0 & \boxed{1} & 0 & 0 & 0 & \\
            \boxed{1} & -1 & 0 & 1 & -1 & 1 & -1 & 1 & \\
            \vdots & \vdots & \vdots & \vdots & \vdots & \vdots &  \vdots &  \vdots\\
            0 & \boxed{1} & 0 & 0 & 0 & 0 & 0 & 0\\
        \end{pmatrix}   $
    \end{center}
    \caption{Part of an ASM with three $-1$s in the same row where there is a $0$ between two nonzero entries in row $i$ left of center and the $1$ in that column is above row $i$ acts as the $1$ (left) or $2$ (right) of a $12$ pattern with the $1$ above the westmost $-1$.}
    \label{F:3neg1s_gaps}
    \end{figure}
    \item The two $1$s in columns $b$ and $\beta$ form a $12$ pattern and $\beta >b$.  Then the $1$ above the $0$ in row $i$ (in column $\beta$),  the $1$ above the center $-1$ (in column $d$), the westmost $1$ of row $i$ (in column $a$), and the $1$ below the westmost $-1$ (in column $b$) form a $3412$ pattern as shown in Figure~\ref{F:3neg1s_gaps} (right).
      \end{enumerate}

    Therefore, there cannot be any $0$ entries between nonzero entries in row $i$.
    \end{proof}

We end the section with the following theorem enumerating ASMs classically avoiding $2143$ and $3412$ that have exactly three negative ones. Recall  $\ASM_{n,k}(\sigma)$ denotes the set of ASMs in $\ASM_n(\sigma)$ containing exactly $k$ negative ones. 
    \begin{ThreeNegativeOnesThm}
        For $n\geq 7$, we have
        \[ |\ASM_{n,3}(2143, 3412)| = \sum_{k=0}^{n-5} 2^{n-4-k} \binom{2k}{k} - (n-2) 2^{n-5}. \]
        For $n \le 6$, we have $|\ASM_{n,3}(2143, 3412)| = 0$. The generating function for this sequence is
        \[ \sum_{n\ge0} |\ASM_{n,3}(2143, 3412)|\,z^n =  \frac{2z^5}{(1-2z)\sqrt{1-4z}} -\frac{2z^5(1-2z^2)}{(1-2z)^2}. \]
    \end{ThreeNegativeOnesThm}

    \begin{proof}
We show that the number of ASMs in $\ASM_n(2143,3412)$ with exactly three negative ones is given by
\[ 2 \sum_{i=4}^{n-3}\sum_{a=1}^{n-6}\sum_{j=1}^{a}\binom{i-2}{j}\binom{n-i-1}{a+1-j}\binom{a-1}{j-1}\binom{n-a-6}{i-j-3}.\]
We prove in Lemma \ref{lem:appendixlemma} in the Appendix that 
the generating function for these terms is 
\[ \frac{2z^5}{(1-2z)\sqrt{1-4z}} -\frac{2z^5(1-2z^2)}{(1-2z)^2}. \]
It is then routine to verify (and we leave this to the reader) that the coefficients of this generating function can be expressed as 
\[ \sum_{k=0}^{n-5} 2^{n-4-k} \binom{2k}{k} - (n-2) 2^{n-5} \]
for $n\ge 7$ and $0$ for $n \le 6$, as desired. 

First, note that by Lemma~\ref{lem:allinarow}, any ASM in $\ASM_n(2143,3412)$ with three $-1$s must have them all in the same row or column. Let $A\in \ASM_n(2143,3412)$ with three $-1$s all in row $i$. There must be a $1$ before the first $-1$, after the last $-1$, and between each $-1$, so $A$ must have size $n \geq 7$. 

    Furthermore, by Lemma~\ref{lem:abcdefg_adjacent}, the columns of all nonzero entries in row $i$ are adjacent. Let $a$ be the column index of the westmost nonzero entry in row $i$. Then, by Lemma~\ref{lem:iplusminus1}, the $1$s above and below the $-1$ in column $a+3$ must be in rows $i$ and $i+1$.  Thus, we know the exact structure of rows $i-1,i,i+1$ of $A$  once we know the value of $a$, and each of the $n-3$ remaining rows have exactly one $1$ in them. The $1$s above row $i-1$ to the west of column $a+3$ (the column with the middle $-1$) must be in increasing order to avoid forming a $21$ pattern that could combined with the eastmost $1$ of row $i$ and the $1$ in row $i+1$ to form a $2143$ pattern.  Similarly, the $1$s above row $i-1$ to the east of column $a+3$ must be in decreasing order to avoid forming a $34$ pattern that can be combined with the westmost $1$ of row $i$ and the $1$ in row $i+1$ to form a $3412$ pattern.  Likewise, the $1$s below row $i+1$ to the west of column $a+3$ must be in decreasing order and the $1$s below row $i+1$ to the east of column $a+3$ must be in increasing order.  
    
    Now form a $0-1$ matrix $M$ consisting of the $n-3$ rows $1,2,\ldots,i-2$ and $i+2,i+3,\ldots,n$ and the $n-5$ columns of the original ASM except for the columns that have a $1$ in row $i$ (columns $a,a+2,a+4,a+6$) and the column that has the center $-1$ (column $a+3$). We consider this matrix as split into quadrants, where the western hemisphere is columns $1,2,\ldots,a-1,a+1$ of $A$  and the northern hemisphere is the first $i-2$ rows. 
    
    Matrix $M$ is similar to a skew-merged permutation matrix~\cite{Atkinson, Stankova} mentioned earlier in this section. 
    The differences are that $M$ specifically has at least one $1$ in each leg of the $X$ shape, every leg is explicitly defined as being in one of the four quadrants with no possible overlap, and columns $a$ and $a+1$ of $M$ (the columns that contained the outer $-1$s  in $A$) each have a $1$ in the northern and southern hemisphere.
    
    We have a total of $i-2$ entries in $M$ that are $1$s in the northern hemisphere and so $n-i-1$ entries that are $1$s in the southern hemisphere.  Since we have two $1$s in column $a$ of $M$, we have a total of $a+1$ entries that are $1$s in the western hemisphere and so there are $n-a-4$ entries that are $1$s in the eastern hemisphere.  Say that there are $j$ entries that are $1$s in the northwest quadrant.  There are then $i-2-j$ entries that are $1$s in the northeast quadrant, $a+1-j$ entries that are $1$s in the southwest quadrant, and $n+j-i-a-2$ entries that are $1$s in the southeast quadrant. This is summarized by the following picture. 

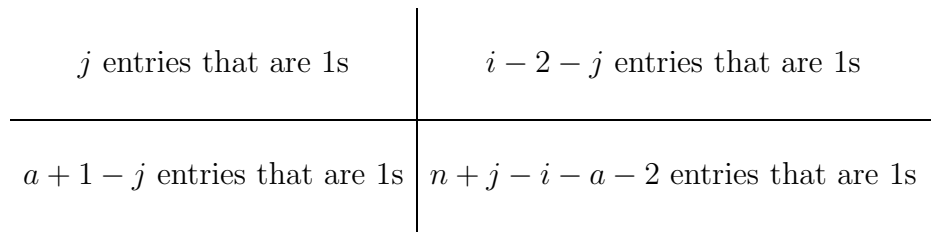
\begin{center}
    \begin{tabular}{c|c}
    &\\
    $j$ entries that are $1$s & $i-2-j$ entries that are $1$s\\
    &\\
    \hline 
    &\\
    $a+1-j$ entries that are $1$s & $n+j-i-a-2$ entries that are $1$s\\
    &\\
    \end{tabular}
    \end{center}

    We must have $n-3$ entries that are $1$s and at least one entry that is a $1$ in each quadrant, so the following must hold:
    \begin{align*}
    2 &\leq i-2 \leq (n-3) -2 \\
    2 & \leq a+1  \leq (n-3)-2 \\
    1 & \leq j  \leq a \\
    1 & \leq i-2-j  \\ 
    1 & \leq n+j-i-a -2. 
    \end{align*}

    Of the $i-2$ rows (which each have one $1$ entry) in the northern hemisphere, there are $\binom{i-2}{j}$ ways to choose which have $1$s in the northwest quadrant.  Similarly, of the $n-i-1$ rows in the southern hemisphere, there are $\binom{n-i-1}{a+1-j}$ ways to choose which have $1$s in the southwest quadrant. 
    Then of the $a+1$ $1$s in the western hemisphere, there is one entry from each quadrant in column $a+1$,
    so we have $\binom{a-1}{j-1}$ ways to select which columns the rest of the $1$ entries are in the northwest quadrant.  Finally, after placing the centermost entry from each of the northeast and southeast quadrants in column $a+1$ of $M$ 
    there are $\binom{n-a-6}{i-j-3}$ ways to select which columns the rest of the $1$s of the northeast quadrant are in.

    To satisfy all the necessary bounds, we can say there are 
    \[\sum_{i=4}^{n-3}\sum_{a=1}^{n-6}\sum_{j=1}^{\min(a,i-3)}\binom{i-2}{j}\binom{n-i-1}{a+1-j}\binom{a-1}{j-1}\binom{n-a-6}{i-j-3} \mathbf{1}[a \leq n+j-i-3]\] 
    ASMs in $\ASM_n(2143,3412)$ with three $-1$s in the same row, where $\mathbf{1}[a \leq n+j-i-3]$ equals $1$ if  $a \leq n+j-i-3$ and $0$ otherwise.

    However, we note that the bounds $j \leq i-3$ and $a \leq n+j-i-3$ are unnecessary, because the last binomial factor $\binom{n-a-6}{i-j-3}$  equals zero if they are not satisfied.  Hence we can also say there are 
    \[\sum_{i=4}^{n-3}\sum_{a=1}^{n-6}\sum_{j=1}^{a}\binom{i-2}{j}\binom{n-i-1}{a+1-j}\binom{a-1}{j-1}\binom{n-a-6}{i-j-3}\] 
    ASMs in $\ASM_n(2143,3412)$ with three $-1$s in the same row.

    Symmetrically, there will be an equal number of ASMs in $\ASM_n(2143,3412)$ with three $-1$s in the same column, and thus the total number of ASMs in $\ASM_n(2143,3412)$ with three $-1$s is 
    $$2 \sum_{i=4}^{n-3}\sum_{a=1}^{n-6}\sum_{j=1}^{a}\binom{i-2}{j}\binom{n-i-1}{a+1-j}\binom{a-1}{j-1}\binom{n-a-6}{i-j-3}.$$
In Lemma \ref{lem:appendixlemma} in the Appendix, we prove that this expression equals the desired formula. 
More precisely, we show that the generating function for these terms is the one announced in the theorem. To the effect, we first ``decouple'' the two instances of $n$ in the above summation, turning our attention (after reindexing) to the bivariate generating function
\[ F(x,y) = \sum_{\gamma = 0}^\infty \sum_{\beta=0}^\infty \sum_{i=0}^\infty \sum_{a=0}^\infty \sum_{j=0}^\infty \binom{i+2}{j+1} \binom{\gamma-i+2}{a-j+1} \binom{a}{j} \binom{\beta-a}{i-j} x^\gamma y^\beta.\] 
Using Wilf's snake oil method and the simple identity $\sum_{n=0}^\infty \binom{n}{k} x^n = \frac{x^k}{(1-x)^{k+1}}$, we find a closed-form expression of the generating function $F(x,y)$. Considering then the rescaling $H(x,z) = F(x,z/x)$, we want to find an expression for the coefficient of $x^0$ in $H(x,z)$.
We express this coefficient as a contour integral, which we evaluate by the method of residues. 
\end{proof}

The first $12$ terms for the sequence of ASMs of size $n\geq 7$ classically avoiding $2143$ and $3412$ with exactly three $-1$s are:
\[8, 48, 220, 912, 3608, 13952, 53388, 203504, 775496, 2959808, 11323832, 43440672.\]
Note that neither this sequence nor the sequence of each term divided by $4$ appeared in the OEIS~\cite{OEIS} at the time of this writing.
 
\section{Future directions}
\label{sec:future}
We consider below some ideas for future research.

\subsection{Other sets of multiple patterns}
It may be interesting to study $\ASM_n(R)$ where the set $R$ consists of one permutation from $S(123)$ and one permutation from $S(321)$. 
For instances, $\ASM_n(321, 2341)$ may be of interest, since the enumeration of $S(321, 2341)$ has a nice formula: $|S_n(321, 2341)| = f_{2n-1}$ (every other Fibonacci number).  Based on our code in Sage~\cite{Sage}, the first few terms of the sequence for $|\ASM_n(321, 2341)|$ are $1, 2, 6, 22, 87, 353, 1445$.
Another interesting open problem would be to study the exact or asymptotic enumeration of $\ASM_n(2413, 3142)$, since $S_n(2413, 3142)$ is the set of separable permutations. 

\subsection{Words and other $0 - 1$ matrices}
Pattern avoidance in ASMs need not be limited to permutations.  One can consider patterns that are \emph{words}, where the corresponding word matrix would still have the rows representing positions and columns giving values.  For example, the word $w=1231$  would be represented as the matrix $M_w=\left[\begin{smallmatrix} 
1 & 0 & 0   \\
0 & 1 & 0   \\
0 & 0 & 1 \\
1 & 0 & 0 
\end{smallmatrix}\right]$. 
In fact, we can consider any $0-1$ matrix as a pattern that can either be classically contained or avoided in the same manner that we considered for permutation matrices.  We include a few introductory cases below.

\begin{prop}
    The ASMs that avoid the word $11$ are exactly the permutation matrices, and so are enumerated by $n!$.
\end{prop}

\begin{proof}
    Two $1$s in the same column gives an instance of a $11$ pattern (word matrix). An ASM contains two $1$s in the same column if and only if there is a $-1$ in that column. The ASMs that do not have a $-1$ are exactly the permutation matrices.
\end{proof}

Using avoidance of words, 
we can describe the ASMs that have at most one $-1$ in terms of pattern avoidance.  Specifically if an ASM has two or more $-1$s, these negative ones could appear in:
\begin{enumerate}
    \item the same column, creating a $111$ pattern,
    \item the same row, creating a $111$ pattern in the transposed matrix, or
    \item in different rows and columns, creating (at least one of) the following patterns: $$1122,1212,1221,2112,2121,2211.$$
\end{enumerate} 

\begin{prop}~\label{p:at_most_one}
    The number of ASMs with at most one $-1$ are exactly the ASMs that avoid all of     $111,1122,1212,1221,2112,2121,2211$, and whose transpose also avoids $111$.  The number of such ASMs of size $n$ is $\frac{n!}{6}\binom{n}{3}+n!$
\end{prop}


\begin{remark}
\label{remark:onenegative}
Note that the first few terms of the sequence of ASMs with \emph{exactly} one $-1$ are $0, 0, 1, 16, 200$, $2400, 29400$, which is sequence A001810 in \cite{OEIS}.
In \cite{Lalonde}, Lalonde found a generating function recurrence for the $q$-enumeration of ASMs with exactly one $-1$, where the power of $q$ corresponds to the inversion number. An explicit formula for the number of ASMs of size $n$ with exactly one $-1$ is given by $\frac{n!}{6}\binom{n}{3}$; see \cite[Appendix A]{AyyerBehrendFischer}. 
\end{remark}
The enumeration given in Proposition~\ref{p:at_most_one} is an immediate consequence as the only other ASMs with at most one $-1$ are the permutation matrices.

\subsection{Avoidance of ASMs and identical avoidance}
One can also consider the avoidance of ASMs rather than permutation matrices, as suggested in \cite[Sec.~6]{JL}.  Denote by $M$ the smallest ASM that is not a permutation matrix, namely $M=\left[\begin{smallmatrix} 
0 & 1 & 0  \\
1 & -1 & 1  \\
0 & 1 & 0 
\end{smallmatrix}\right]$.

A choice consistent with the conventions of this paper is to consider an ASM $A$ to avoid an ASM $B$ if the $1$s in $A$ avoid the pattern of $1$s in $B$. That is, we may treat each  $-1$ in the pattern matrix as a $0$ and consider avoidance of this $0 - 1$ matrix in the sense of the prior subsection. For instance, $A$ avoids $M$ whenever $A$ avoids $\left[\begin{smallmatrix} 
0 & 1 & 0  \\
1 & 0 & 1  \\
0 & 1 & 0 
\end{smallmatrix}\right]$.
In this case, the ASMs avoiding $M$ are exactly those ASMs which are also permutation matrices.   Furthermore, any ASM that is not a permutation matrix is avoided by all permutation matrices.
\begin{prop}
    For all natural numbers $n$, we have $| \ASM_n(M) | = n!$.  Furthermore,  $| \ASM_n(A) | \geq n!$, where $A$ is any ASM containing at least one $-1$.
\end{prop} 

However, it may make more sense to consider \emph{identical avoidance} where the submatrix must match the pattern ASM exactly.  Let $\ASM^{\text{iden}}_n(R)$ (resp.\ $\ASM^{\text{iden}}(R)$) denote the set of ASMs in $\ASM_n$ (resp.\ $\ASM$) avoiding every element of $R$ in this exact (or identical) sense.   Using Sage~\cite{Sage}, we have $$|\ASM_n^{\text{iden}}(M)| = 1,2,6,26,176,1886, 29088,\ldots$$  

Either notion of avoidance could similarly be applied to ASMs or even to more general matrices on the same alphabet of $\{-1,0,1\}$.

\begin{remark}  As identical containment is a stricter definition, there are at least as many ASMs of size $n$ identically avoiding any ASM $A$ (or word $w$) than ASMs of the same size which classically avoid $A$ (respectively $w$).
\end{remark} 

\section*{Acknowledgments}
This work began at the 2023 Schloss Dagstuhl workshop on Pattern Avoidance, Statistical Mechanics and Computational Complexity. We wish to thank the institute and the workshop organizers for the conducive research environment. The work was later supported by opportunities to present and collaborate at the Permutation Patterns Conferences in both 2023 and 2024. The authors also thank Ilse Fischer for pointing us to reference \cite{AyyerBehrendFischer} for Remark~\ref{remark:onenegative} and Alex Woo for helpful conversations. Finally, we are grateful to the reviewers for their constructive feedback on our work.  Striker was supported by a Simons Foundation gift MP-TSM-00002802 and NSF grant DMS-2247089.

\section*{Appendix: Binomial coefficient identity for Theorem \ref{thm:exactly3enumeration}}\label{secA1}

To conclude the proof of Theorem \ref{thm:exactly3enumeration}, we use an identity to simplify the formula we obtained.
Our proof of this identity is computationally intensive, so we have included it here rather than in the proof of the theorem.

\begin{lem}
\label{lem:appendixlemma}
For all $n \ge 7$ we have
\[2 \sum_{i=4}^{n-3} \sum_{a=1}^{n-6} \sum_{j=1}^a \binom{i-2}{j} \binom{n-i-1}{a+1-j} \binom{a-1}{j-1} \binom{n-a-6}{i-j-3} = \sum_{k=0}^{n-5} 2^{n-4-k} \binom{2k}{k} - (n-2) 2^{n-5}. \]
\end{lem}
\begin{proof}
    If we reindex the sums on the left and divide by two, then we see that it is sufficient to show that
    \begin{equation}
        \label{eqn:reindexed}
    \sum_{i=0}^n \sum_{a=0}^n \sum_{j=0}^a \binom{i+2}{j+1} \binom{n-i+2}{a-j+1} \binom{a}{j} \binom{n-a}{i-j} = \sum_{k=0}^{n+2} 2^{n+2-k} \binom{2k}{k} - (n+5) 2^{n+1}.
    \end{equation}
    This is equivalent to showing that the generating function for the sequence on the left is 
    \[ \frac{1}{z^2 (1-2z)\sqrt{1-4z}} - \frac{1-2z^2}{z^2(1-2z)^2}, \]
    since this is the generating function for the sequence on the right.

    To begin, consider the generating function in $x$ and $y$ given by
    \[ F(x,y) = \sum_{\gamma = 0}^\infty \sum_{\beta=0}^\infty \sum_{i=0}^\infty \sum_{a=0}^\infty \sum_{j=0}^\infty \binom{i+2}{j+1} \binom{\gamma-i+2}{a-j+1} \binom{a}{j} \binom{\beta-a}{i-j} x^\gamma y^\beta.\]
    Note that if we adopt the convention that $\binom{n}{k} \neq 0$ implies $n \ge k \ge 0$ then the coefficients of $x^n y^n$ in $F(x,y)$ are the terms of the sequence on the left side of \eqref{eqn:reindexed}.
    We will use Wilf's snake oil method \cite[Sec.~4.3]{Wilf} to obtain a closed form expression for $F(x,y)$.
    This will be a rational function of $x$ and $y$, and we will then be able to use techniques outlined in \cite[Sec.~6.3]{StanleyVol2} to obtain the desired generating function.

    In the computations that follow we will use the fact that if $k$ is a nonnegative integer then
\begin{equation}
\label{eqn:thesum}
\sum_{n=0}^\infty \binom{n}{k} x^n = \frac{x^k}{(1-x)^{k+1}}.
\end{equation}

Reordering the sums in our expression for $F$, we have
\[ F(x,y) = \sum_{\gamma = 0}^\infty \sum_{i=0}^\infty \sum_{a=0}^\infty \sum_{j=0}^\infty \binom{i+2}{j+1}  \binom{\gamma-i+2}{a-j+1} \binom{a}{j} x^\gamma y^a \sum_{\beta = 0}^\infty \binom{\beta-a}{i-j} y^{\beta-a}. \]
The rightmost sum is zero if $i < j$, so we should start the sum over $i$ with $i-j$.
Now when we apply \eqref{eqn:thesum} this becomes
\begin{eqnarray*}
F(x,y) &=& \sum_{\gamma = 0}^\infty \sum_{j=0}^\infty \sum_{a=0}^\infty \sum_{i=j}^\infty \binom{i+2}{j+1}  \binom{\gamma-i+2}{a-j+1} \binom{a}{j} x^\gamma y^a \sum_{\beta = 0}^\infty \binom{\beta-a}{i-j} y^{\beta-a} \\
&=& \sum_{\gamma = 0}^\infty \sum_{j=0}^\infty \sum_{a=0}^\infty \sum_{i=j}^\infty \binom{i+2}{j+1}  \binom{\gamma-i+2}{a-j+1} \binom{a}{j} x^\gamma y^a \frac{y^{i-j}}{(1-y)^{i-j+1}}.
\end{eqnarray*}
Reindexing the innermost sum with $k = i-j$, so $i = k + j$, we have
\begin{eqnarray*}
F(x,y) &=& \sum_{\gamma = 0}^\infty \sum_{j=0}^\infty \sum_{a=0}^\infty \sum_{k=0}^\infty \binom{k+j+2}{j+1}  \binom{\gamma-k-j+2}{a-j+1} \binom{a}{j} x^\gamma y^a \frac{y^k}{(1-y)^{k+1}} \\
&=&  \sum_{j=0}^\infty \sum_{a=0}^\infty \sum_{k=0}^\infty \binom{k+j+2}{j+1}  \binom{a}{j} y^a \frac{y^k}{(1-y)^{k+1}}   \sum_{\gamma = 0}^\infty  \binom{\gamma-k-j+2}{a-j+1} x^\gamma.
\end{eqnarray*}
The rightmost sum is zero if $a < j-1$ and $\binom{a}{j} = 0$ if $a = j-1$, so this becomes
\[ F(x,y) = \sum_{j=0}^\infty \sum_{a=j}^\infty \sum_{k=0}^\infty \binom{k+j+2}{j+1}  \binom{a}{j} \frac{y^{a+k}}{(1-y)^{k+1}} x^{k+j-2} \sum_{\gamma = 0}^\infty  \binom{\gamma-k-j+2}{a-j+1} x^{\gamma-k-j+2}. \]
Reindexing the innermost sum with $\omega = \gamma - k - j + 2$, so $\gamma = \omega + j + k -2$, we find
\[ F(x,y) = \sum_{j=0}^\infty \sum_{a=j}^\infty \sum_{k=0}^\infty \binom{k+j+2}{j+1}  \binom{a}{j} \frac{y^{a+k}}{(1-y)^{k+1}} x^{k+j-2} \sum_{\omega = 2-k-j}^\infty  \binom{\omega}{a-j+1} x^\omega. \]

If $k+j \le 1$ then the inner sum begins with $\omega > 0$, and \eqref{eqn:thesum} does not apply directly.
Now if $k+j = 1$ then the inner sum begins with $\omega =1$.
In this case the only term that can be missing occurs when $a -j +1 = 0$, which happens when $a = j-1$.
But $a \ge j$, so this term is zero.
With all of this in mind, we separate the terms in our sum into four cases.

\medskip

\noindent
{\bf Case One:}  $j = k = 0$.

\smallskip

\noindent
The generating function for such terms is
\begin{align*}
\sum_{a=0}^\infty 2 \frac{y^a}{x^2 (1-y)} & \sum_{\omega = 2}^\infty  \binom{\omega}{a+1} x^{\omega} \\
&= \frac{2}{x^2(1-y)} \sum_{\omega=2}^\infty \binom{\omega}{1} x^{\omega} + \sum_{a=1}^\infty 2 \frac{y^a}{x^2(1-y)} \sum_{\omega=2}^\infty \binom{\omega}{a+1} x^{\omega} \\
&= \frac{2}{x^2(1-y)} \left( \frac{x}{(1-x)^2} - x \right) + \sum_{a=1}^\infty 2 \frac{y^a}{x^2(1-y)} \frac{x^{a+1}}{(1-x)^{a+2}} \\
&= \frac{2}{x^2(1-y)} \left( \frac{x}{(1-x)^2} - x \right) + \frac{2}{x(1-y)(1-x)^2} \sum_{a=1}^\infty \left(\frac{xy}{1-x}\right)^a \\
&= \frac{2}{x^2(1-y)} \left( \frac{x}{(1-x)^2} - x \right) + \frac{2y}{(1-y)(1-x)^2(1-x-xy)} \\
&= \frac{2}{x(1-y)(1-x)^2} - \frac{2}{x(1-y)} + \frac{2y}{(1-y)(1-x)^2(1-x-xy)}.
\end{align*}

\medskip

\noindent
{\bf Case Two:} $j =0$ and $k \ge 1$.

\smallskip

\noindent
The generating function for such terms is
\begin{align*}
\sum_{a=0}^\infty \sum_{k=1}^\infty & (k+2) \frac{y^{a+k}}{(1-y)^{k+1}} x^{k-2} \sum_{\omega = 2-k}^\infty \binom{\omega}{a+1} x^\omega \\
&= \sum_{a=0}^\infty \sum_{k=1}^\infty (k+2) \frac{y^{a+k}}{(1-y)^{k+1}} x^{k-2} \frac{x^{a+1}}{(1-x)^{a+2}} \\
&= \left(\sum_{a=0}^\infty \frac{y^a x^{a+1}}{(1-x)^{a+2}} \right) \left( \sum_{k=1}^\infty (k+2) \frac{y^k x^{k-2}}{(1-y)^{k+1}} \right) \\
&= \frac{1}{x (1-x)^2 (1-y)} \left(\sum_{a=0}^\infty \left(\frac{xy}{1-x}\right)^a \right) \left( \sum_{k=1}^\infty (k+2) \left(\frac{xy}{1-y}\right)^k \right) \\
&= \frac{1}{x (1-x)^2 (1-y)} \left( \frac{1}{1-\frac{xy}{1-x}} \right) \left(\frac{1-y}{xy} \right)^2 \left( \sum_{k=1}^\infty \binom{k+2}{1} \left(\frac{xy}{1-y}\right)^{k+2} \right) \\
&= \frac{1-y}{x^3 y^2 (1-x)(1-x-xy)} \left( \frac{\frac{xy}{1-y}}{\left(1-\frac{xy}{1-y}\right)^2} - \frac{xy}{1-y} -2\left(\frac{xy}{1-y}\right)^2\right)\\
&= \frac{1-y}{x^3 y^2 (1-x)(1-x-xy)} \left(\frac{xy(1-y)}{(1-y-xy)^2} - \frac{xy}{1-y} -2\left(\frac{xy}{1-y}\right)^2\right) \\
&= \frac{(1-y)^2}{x^2 y (1-x)(1-x-xy)(1-y-xy)^2} - \frac{1}{x^2 y (1-x)(1-x-xy)} \\
& \hspace{150pt} - \frac{2}{x (1-x) (1-y)(1-x-xy)}.
\end{align*}

\medskip

\noindent
{\bf Case Three:}. $j \ge 1$ and $k = 0$.

\smallskip

\noindent
The generating function for such terms is
\begin{align*}
\sum_{j=1}^\infty & \sum_{a=j}^\infty (j+2) \binom{a}{j} \frac{y^a}{1-y} x^{j-2} \sum_{\omega=2-j}^\infty \binom{\omega}{a-j+1} x^\omega \\
&= \sum_{j=1}^\infty \sum_{a=j}^\infty (j+2) \binom{a}{j} \frac{y^a}{1-y} x^{j-2} \left( \frac{x^{a-j+1}}{(1-x)^{a-j+2}} \right) \\
&= \sum_{j=1}^\infty (j+2) \frac{(1-x)^{j-2}}{x(1-y)} \sum_{a=j}^\infty \binom{a}{j} \left( \frac{x y}{1-x} \right)^a \\
&= \sum_{j=1}^\infty (j+2) \frac{(1-x)^{j-2}}{x(1-y)} \left( \frac{(1-x) (xy)^j}{(1-x-xy)^{j+1}}\right)\\
&= \frac{1}{x(1-x)(1-y)(1-x-xy)} \sum_{j=1}^\infty (j+2) \left( \frac{xy(1-x)}{1-x-xy} \right)^j\\
&= \frac{1}{x(1-x)(1-y)(1-x-xy)} \left( \frac{1-x-xy}{xy(1-x)}\right)^2 \sum_{j=1}^\infty \binom{j+2}{1} \left( \frac{xy(1-x)}{1-x-xy} \right)^{j+2}\\
&= \frac{1-x-xy}{x^3 y^2(1-x)^3(1-y)} \left( \frac{\frac{xy(1-x)}{1-x-xy}}{\left(1-\frac{xy(1-x)}{1-x-xy}\right)^2} - \frac{xy(1-x)}{1-x-xy} -2 \left( \frac{xy(1-x)}{1-x-xy}\right)^2\right) \\
&= \frac{1-x-xy}{x^3 y^2(1-x)^3(1-y)} \left( \frac{xy(1-x)(1-x-xy)}{(1-x-2xy+x^2y)^2}- \frac{xy(1-x)}{1-x-xy} -2 \left( \frac{xy(1-x)}{1-x-xy}\right)^2\right) \\
&= \frac{(1-x-xy)^2}{x^2 y(1-x)^2 (1-y)(1-x-2xy+x^2y)^2} - \frac{1}{x^2y (1-x)^2 (1-y)} \\
& \hspace{170pt} - \frac{2}{x (1-x)(1-y)(1-x-xy)}.
\end{align*}

\medskip

\noindent
{\bf Case Four:}
$j \ge 1$ and $k \ge 1$.

\smallskip

\noindent
The generating function for such terms is
\begin{align*}
 \sum_{j=1}^\infty \sum_{a=j}^\infty & \sum_{k=1}^\infty \binom{k+j+2}{j+1} \binom{a}{j} \frac{y^{a+k}}{(1-y)^{k+1}} x^{k+j-2} \sum_{\omega = 2-k-j}^\infty  \binom{\omega}{a-j+1} x^\omega \\
 &=  \sum_{j=1}^\infty \sum_{a=j}^\infty \sum_{k=1}^\infty \binom{k+j+2}{j+1} \binom{a}{j} \frac{y^{a+k}}{(1-y)^{k+1}} x^{k+j-2} \left( \frac{x^{a-j+1}}{(1-x)^{a-j+2}} \right) \\
 &=  \sum_{j=1}^\infty \sum_{k=1}^\infty \binom{k+j+2}{j+1}\frac{y^{k}}{(1-y)^{k+1}}  \left( \frac{x^{k-1}}{(1-x)^{-j+2}} \right)  \sum_{a=j}^\infty \binom{a}{j}\left( \frac{x y}{1-x} \right)^a \\
 &=  \sum_{j=1}^\infty \sum_{k=1}^\infty \binom{k+j+2}{j+1}\frac{y^{k}}{(1-y)^{k+1}}  \left( \frac{x^{k-1}}{(1-x)^{-j+2}} \right) \frac{x^j y^j (1-x)}{(1-x-xy)^{j+1}} \\
 &= \sum_{j=1}^\infty \frac{(1-x)^{j-1} x^{j-1} y^j}{(1-x-xy)^{j+1} (1-y)} \sum_{k=1}^\infty \binom{k+j+2}{j+1} \left( \frac{xy}{1-y} \right)^k \\
 &= \sum_{j=1}^\infty \frac{(1-x)^{j-1} (1-y)^{j+1}}{x^3 y^2 (1-x-xy)^{j+1} } \sum_{k=1}^\infty \binom{k+j+2}{j+1} \left( \frac{xy}{1-y} \right)^{k+j+2} \\
 &= \sum_{j=1}^\infty \frac{(1-x)^{j-1} (1-y)^{j+1}}{x^3 y^2 (1-x-xy)^{j+1} } \left( \frac{x^{j+1} y^{j+1} (1-y)}{(1-y-xy)^{j+2}} - \left(\frac{xy}{1-y} \right)^{j+1} - (j+2) \left(\frac{xy}{1-y}\right)^{j+2} \right) \\
 &= \sum_{j=1}^\infty \frac{(1-x)^{j-1} (1-y)^{j+2} x^{j-2} y^{j-1}}{(1-x-xy)^{j+1} (1-y-xy)^{j+2}} - \sum_{j=1}^\infty \frac{(1-x)^{j-1} x^{j-2} y^{j-1}}{(1-x-xy)^{j+1}} \\
 & \hspace{160pt} - \sum_{j=1}^\infty \frac{(j+2) (1-x)^{j-1} x^{j-1} y^j}{(1-x-xy)^{j+1}(1-y)}.
 \end{align*}
 
The leftmost sum is
\begin{align*}
\frac{(1-y)^3}{x (1-x-xy)^2 (1-y-xy)^3} & \sum_{j=1}^\infty \left(\frac{(1-x)(1-y) x y}{(1-x-xy)(1-y-xy)}\right)^{j-1}  \\
&= \frac{(1-y)^3}{x (1-x-xy)^2 (1-y-xy)^3}  \left( \frac{1}{1- \frac{(1-x)(1-y) x y}{(1-x-xy)(1-y-xy)}} \right) \\
&= \frac{(1-y)^3}{x (1-x-xy) (1-y-xy)^2 (1-x-y-2xy+2x^2y+2xy^2)}
\end{align*}

The middle sum is
\begin{align*}
\frac{1}{x(1-x-xy)^2} & \sum_{j=1}^\infty \left(\frac{(1-x)xy}{1-x-xy}\right)^{j-1} \\
&= \frac{1}{x(1-x-xy)^2}\left( \frac{1}{1-\frac{(1-x)xy}{1-x-xy}} \right) \\
&= \frac{1}{x(1-x-xy) (1-x-2xy+x^2y)}\\
\end{align*}

The rightmost sum is
\begin{align*}
& \frac{1-x-xy}{x^3 y^2 (1-x)^3(1-y)} \sum_{j=1}^\infty \binom{j+2}{1} \left( \frac{xy(1-x)}{1-x-xy}\right)^{j+2} \\
&= \frac{1-x-xy}{x^3 y^2 (1-x)^3(1-y)}\left(\frac{xy(1-x)(1-x-xy)}{(1-x-2xy+x^2y)^2} - \frac{xy(1-x)}{1-x-xy} -2 \left(\frac{xy(1-x)}{1-x-xy}\right)^2 \right) \\
&= \frac{(1-x-xy)^2}{x^2y(1-x)^2(1-y)(1-x-2xy+x^2y)^2} -\frac{1}{x^2y(1-x)^2(1-y)} \\
& \hspace{140pt} -\frac{2}{x(1-x)(1-y)(1-x-xy)}
\end{align*}

When we combine all of these cases we find
\[ F(x,y) = \frac{G(x,y)}{(1-x)(1-y)(1-y-xy)(1-x-2xy+x^2y)(1-x-y)(1-2xy)}, \]
where
\begin{align*}
G(x,y) = 4-6x+ & 2x^2-18xy+24x^2y-8x^3y-9y^2+20xy^2+14x^2y^2-30x^3y^2 \\
& +10x^4y^2 +5y^3+10xy^3-32x^2y^3+6x^3y^3+12x^4y^3-4x^5y^3 \\
& -8xy^4+4x^2y^4+8x^3y^4-4x^4y^4.
\end{align*}

To obtain the generating function for this sequence, we rescale $F$ by setting
\[ H(x,z) = F\left(x,\frac{z}{x}\right). \]
Now $H$ is a Laurent series in $x$ and the generating function we want is the ``constant'' term in $H$, which is the coefficient of $x^0$.

To begin, we have
\[ H(x,z) = \frac{J(x,z)}{(1-x)(z-x)(1-x-2z+xz)(z-x+x^2)(x-xz-z)(1-2z)} \]
where
\begin{align*}
J(x,z) = 4x^3 &- 6x^4 +2x^5-18x^3z+24x^4z-8x^5z-9xz^2+20x^2z^2+14x^3z^2 \\
& -30x^4z^2+ 10x^5z^2+5z^3+10xz^3-32x^2z^3+6x^3z^3+12x^4z^3-4x^5z^3 \\
& -8z^4+4xz^4+8x^2z^4-4x^3z^4.
\end{align*}
We now use the fact that our desired generating function can be expressed as the contour integral
\begin{align*}
    \frac{1}{2\pi i} \int_C & \frac{H(x,z)}{x}\; dx \\
    & = \frac{1}{2\pi i} \int_C \frac{J(x,z)}{x(1-x)(z-x)(1-x-2z+xz)(z-x+x^2)(x-xz-z)(1-2z)}\; dx.
\end{align*}
Here $C$ is a circle centered at the origin with radius less than one, and we take $z$ to be near zero to ensure all of the series involved converge.

The poles for the integrand are $0$, $1$, $z$, $\frac{1-2z}{1-z}$, $\frac{1-\sqrt{1-4z}}{2}$, $\frac{1+\sqrt{1-4z}}{2}$ and $\frac{z}{1-z}$, and all of these poles are simple.
Since the radius of $C$ is less than one and $z$ is near zero, the only poles that are inside $C$ are those which go to zero as $z$ goes to zero.
These are $0$, $z$, $\frac{1-\sqrt{1-4z}}{2}$, and $\frac{z}{1-z}$.
The associated residues are, respectively, $\frac{8z-5}{(1-2z)^2}$, $-\frac{2}{z}$, $\frac{1}{z^2(1-2z)\sqrt{1-4z}}$, and $-\frac{(1-z)^2}{z^2(1-2z)^2}$.
The value of the integral, and therefore our generating function, is the sum of these residues.
This is $\frac{1}{z^2 (1-2z)\sqrt{1-4z}} - \frac{1-2z^2}{z^2(1-2z)^2}$, which proves our result.
\end{proof}


\begin{thebibliography}{99}

\bibitem{Atkinson} M. D. Atkinson, Permutations which are the union of an increasing and a decreasing subsequence, \emph{Electron. J. Combin.}, \textbf{5} (1998): \#RP6:1–13.

\bibitem{Arratia}  R. Arratia, On the Stanley-Wilf Conjecture for the Number of Permutations Avoiding a Given Pattern, \emph{Electron. J. Combin.}, \textbf{6} (1999): \#N1:1–4.

\bibitem{AyyerBehrendFischer}
A. Ayyer, R. Behrend, and I. Fischer.
Extreme diagonally and antidiagonally symmetric alternating sign matrices of odd order, \emph{Adv. Math.} \textbf{367} (2020), 107125, 56 pp.

\bibitem{ACG} A. Ayyer, R. Cori, and D. Gouyou-Beauchamps, Monotone triangles and $312$ pattern avoidance, \emph{Electron. J. Combin.} \textbf{18(2)} (2011): \#P26.

\bibitem{BWX} J. Backelin, J. West, and G. Xin, Wilf-equivalence for singleton classes, \emph{Adv. in Appl. Math.} \textbf{38} (2007), no. 2, 133--148.

\bibitem{Bevan} D. Bevan, Permutation patterns: basic definitions and notation, \href{https://arxiv.org/abs/1506.06673}{\texttt{arXiv:1506.06673}} (2015).

\bibitem{Bevan1324} D. Bevan, R. Brignall, A. Price, and J. Pantone, A structural characterisation of $Av(1324)$ and new bounds on its growth rate, \emph{European J. Combin.} \textbf{88} (2020), 103115, 29 pp.

\bibitem{BloomElizalde} J. Bloom and S. Elizalde, Pattern avoidance in matchings and partitions, \emph{Electron. J. Combin.} \textbf{20} (2013): \#P5.

\bibitem{Bona} M. B\'ona, \emph{Combinatorics of Permutations}, Discrete Mathematics and its Applications (Boca Raton), Chapman \& Hall/CRC (2004), Boca Raton, FL.

\bibitem{Bona1342} M. B\'ona, Exact enumeration of $1342$-avoiding permutations: a close link with labeled trees and planar maps, \emph{J. Combin. Theory Ser. A} \textbf{80} (1997), no. 2, 257-272.

\bibitem{Bonairrational} M. B\'ona, The limit of a Stanley-Wilf sequence is not always rational, and layered patterns beat monotone patterns, \emph{J. Combin. Theory Ser. A} \textbf{110} (2005), no. 2, 223--235.

\bibitem{Bonarecords} M. B\'ona, New records in Stanley--Wilf limits, \emph{European J. Combin.} \textbf{28} (2007), no. 1, 75--85.

\bibitem{keyAv} M. Bouvel, R. Smith, and J. Striker, Key-avoidance for alternating sign matrices, 
\emph{Discrete Math. Theor. Comput. Sci.} \textbf{27} (2025), no. 1, 28 pp.

\bibitem{BressoudBook} D. M. Bressoud, \emph{Proofs and Confirmations: The Story of the Alternating-Sign Matrix Conjecture}. Spectrum.
Cambridge University Press, 1999.

\bibitem{BKMS} R. A. Brualdi, K. P. Kiernan, S. A. Meyer, and M. W. Schroeder, Patterns of alternating sign matrices, \emph{Linear Algebra Appl.} \textbf{438} (2013): 3967--3990.

\bibitem{ES} P. Erd\H{o}s and G. Szekeres, A combinatorial problem in geometry. {\em Compositio Math.}. \textbf{2} pp. 463-470 (1935), \url{http://www.numdam.org/item/CM_1935__2__463_0.pdf}.

\bibitem{FK_ASM}
I. Fischer, M. Konvalinka, 
The mysterious story of square ice, piles of cubes, and bijections. \emph{Proc. Natl. Acad. Sci.} (2020), no.38, 23460-23466.

\bibitem{grid} S. Huczynska and V. Vatter, Grid Classes and the Fibonacci Dichotomy for Restricted Permutations, \emph{Electron. J. Combin.}, \textbf{13} (2006): \#R54:1–14.

\bibitem{JL} R. Johansson and S. Linusson, Pattern avoidance in alternating sign matrices, \emph{Ann. Comb.} \textbf{11} (2007): 471--480.

\bibitem{Kitaev} S. Kitaev, \emph{Patterns in Permutations and Words}, Springer, 2011.

\bibitem{Knuth} D. E. Knuth, \emph{The Art Of Computer Programming Vol.\ 1} (Section 2.2.1, Exercises 4 and 5), Addison-Wesley, 1968. 

\bibitem{Lalonde} P. Lalonde, q-Enumeration of alternating sign matrices with exactly one $-1$, \emph{Discrete Math.}
\textbf{256}, no.~3
(2002): 759--773.

\bibitem{MacMahon}
P. A. MacMahon, \emph{Combinatory Analysis}, London: Cambridge University Press, 1916.

\bibitem{MarcusTardos} A. Marcus and G. Tardos, Excluded permutation matrices and the Stanley–Wilf conjecture, {\em J. Combin. Theory Ser. A}. \textbf{107}(1), 153–160 (2004), \url{https://doi.org/10.1016/j.jcta.2004.04.002}

\bibitem{Regev}
A. Regev, Asymptotic values for degrees associated with strips of Young diagrams, {\em Advances in Mathematics}. \textbf{31} (1981), 115--136.

\bibitem{Sage} Sage Mathematics Software, The Sage Development Team (version 10.3). \url{http://www.sagemath.org}.

\bibitem{SS85} 
R. Simion and F. W. Schmidt, Restricted permutations, {\em European Journal Of Combinatorics}, 6:383-406, 1985.

\bibitem{OEIS} 
N. J. A. Sloane, editor, 
\newblock \emph{The On-Line Encyclopedia of Integer Sequences}, published electronically at \url{https://oeis.org}, 2022.

\bibitem{Stankova} Z. E. Stankova, Forbidden subsequences, \emph{Discrete Mathematics}, 132 (1–3) (1994): 291–316, doi:10.1016/0012-365X(94)90242-9.

\bibitem{StanleyVol2} R. P. Stanley, \emph{Enumerative Combinatorics}, Vol.~2, Cambridge University Press, 1999.

\bibitem{Striker} J. Striker, The alternating sign matrix polytope, \emph{Electron. J. Combin.}, \textbf{16} (2009).

\bibitem{Vatter} V. Vatter, Permutation  classes,  \emph{Discrete  Math. Appl. (Boca  Raton)}, CRC Press, Boca Raton, FL, 2015.

\bibitem{Wilf} H. S. Wilf, \emph{generatingfunctionology}, AK Peters, 3rd edition, 2006.

\end{thebibliography}
\end{document}